\documentstyle[12pt]{article}
\font\diagfont=msam10
\makeatletter
\def\carrow{\mbox{\diagfont\char'010}}

\def\LRA#1#2{\@tempdimb=\c@enumiv\@tempdima%
   \vcenter{\offinterlineskip\halign{##\cr%
   \hfil${\scriptstyle{#1}}$\hfil\crcr%
   \hbox to \@tempdimb{\rightarrowfill}\cr%
   \noalign{\kern-1ex}%
   \hbox to \@tempdimb{\leftarrowfill}\cr%
   \hfil${\scriptstyle{#2}}$\hfil\crcr}}}
\def\RA#1{\@tempdimb=\c@enumiv\@tempdima\vbox{\offinterlineskip%
   \halign{##\cr\hfil${\scriptstyle {#1}}$\hfil\crcr%
   \hbox to \@tempdimb{\rightarrowfill}\cr}}}
\def\LA#1{\@tempdimb=\c@enumiv\@tempdima\vbox{\offinterlineskip%
   \halign{##\cr\hfil${\scriptstyle {#1}}$\hfil\crcr%
   \hbox to \@tempdimb{\leftarrowfill}\cr}}}

\def\diag{\leavevmode\bgroup\setcounter{enumiv}{1}%
   \unitlength1em \@tempdima3em \def\\{\crcr&}\vbox\bgroup%
   \def\multicolumn##1##2{\multispan##1\setcounter{enumiv}{##1}%
   \hfil{##2}\hfil\setcounter{enumiv}{1}}
   \offinterlineskip\halign\bgroup\vrule height.8em depth.7em %
   width0pt##&&\hfil${\displaystyle{##}}$\hfil\cr&}
\def\enddiag{\crcr\egroup\egroup\egroup}
\makeatother
\if@twoside \oddsidemargin -20pt \evensidemargin 20pt \marginparwidth 85pt
\else \oddsidemargin 0pt \evensidemargin 0pt
\fi
\advance\textheight by \topskip
\font\symbolfont=msbm10
\font\p=msbm10 at 12pt
\font\teneu=eufm10
\font\egteu=eufm8
\def\dn#1{\mathchoice{\hbox{\teneu #1}}{\hbox{\teneu #1}}%
   {\hbox{\egteu #1}}{\hbox{\egteu #1}}}
\font\symbolfont=msbm10
\def\aaa{\mathop{A}\nolimits}
\def\ddd{\mathop{D}\nolimits}

\def\ccc{\hbox{\symbolfont C}}
\def\fff{\hbox{\symbolfont F}}
\def\nnn{\hbox{\symbolfont N}}
\def\qqq{\hbox{\symbolfont Q}}
\def\rrr{\hbox{\symbolfont R}}

\def\zzz{\hbox{\symbolfont Z}}
\def\cc{\mathop{\cal C}\nolimits}
\def\ff{\mathop{\cal F}\nolimits}

\def\ll{\mathop{\dn{L}}\nolimits}
\def\ma{\mathop{\rm M}\nolimits}
\def\ind{\mathop{\rm ind}\nolimits}
\def\ob{\mathop{\rm ob}\nolimits}
\def\lat{\mathop{\rm lat}\nolimits}
\def\mod{\mathop{\rm mod}\nolimits}
\def\fin{\mathop{\rm fin}\nolimits}
\def\length{\mathop{\rm length}\nolimits}
\def\cm{\mathop{\rm CM}\nolimits}
\def\pr{\mathop{\rm pr}\nolimits}
\def\rin{\mathop{\rm in}\nolimits}
\def\pd{\mathop{\rm pd}\nolimits}
\def\hom{\mathop{\rm Hom}\nolimits}

\def\endm{\mathop{\rm End}\nolimits}

\def\add{\mathop{\rm add}\nolimits}
\def\real{\mathop{\rm Re}\nolimits}
\def\ovo#1{\mathop{\dn{O}_0(#1)}\nolimits}
\def\ovr#1{\mathop{\dn{O}(#1)}\nolimits}
\def\hall{\mathop{\cal H}\nolimits}

\def\hallm{\mathop{\cal M}\nolimits}

\def\plim{\mathop{\raise-.2em\hbox{$\def\arraystretch{0}\begin{array}{c}\lim\\ {\scriptstyle\stackrel{\longleftarrow}{}}\end{array}$}}}
\font\p=msbm10 at 12pt

\def\subsetneq{\mathop{\mbox{ {\p\char'050} }}\nolimits}
\def\Q#1{\begin{picture}(0,0)\put(-1,-2){$#1$}\end{picture}}

\def\QDL{\begin{picture}(7,0)\put(15,5){\vector(-1,-1){10}}\end{picture}}
\def\QDR{\begin{picture}(7,0)\put(5,5){\vector(1,-1){10}}\end{picture}}
\def\QD{\begin{picture}(7,0)\put(5,5){\vector(0,-1){10}}\end{picture}}

\def\ZA{1}
\def\ZAA{1.1}
\def\ZAB{1.2}
\def\ZAC{1.3}
\def\ZAD{1.4}
\def\ZAE{1.5}
\def\ZAF{1.6}
\def\ZAG{1.7}
\def\ZB{2}
\def\ZBA{2.1}
\def\ZBB{2.2}
\def\ZBC{2.3}
\def\ZBCA{2.3.1}
\def\ZBD{2.4}
\def\ZBE{2.5}
\def\ZBF{2.6}
\def\ZBG{2.7}
\def\ZBGA{2.7.1}
\def\ZBGB{2.7.2}
\def\ZC{3}
\def\ZCA{3.1}
\def\ZCB{3.2}
\def\ZCC{3.3}
\def\ZCD{3.4}
\def\ZCDA{3.4.1}
\def\ZCDB{3.4.2}
\def\ZCDC{3.4.3}
\def\ZCE{3.5}
\def\ZCF{3.6}
\def\ZCFA{3.6.1}
\def\ZCFB{3.6.2}
\def\ZCG{3.7}
\def\ZCGA{3.7.1}
\def\ZCGB{3.7.2}
\def\ZD{4}
\def\ZDA{4.1}
\def\ZDB{4.2}
\def\ZDC{4.3}
\def\ZDD{4.4}
\def\ZDE{4.5}

\def\ZDF{4.6}
\def\ZDG{4.7}
\def\ZDGA{4.7.1}
\def\ZDGB{4.7.2}
\def\ZDGC{4.7.3}
\def\ZE{5}
\def\ZEA{5.1}
\def\ZEB{5.2}
\def\ZEBA{5.2.1}
\def\ZEBB{5.2.2}
\def\gl{\mathop{\rm gl.dim}\nolimits}
\def\rdim{\mathop{\rm rep.dim}\nolimits}
\def\fdim{\mathop{\rm fin.dim}\nolimits}
\def\domdim{\mathop{\rm dom.dim}\nolimits}

\def\tri{\mathop{\rm T}\nolimits}
\def\pr{\mathop{\rm pr}\nolimits}
\begin{document}
\begin{center} 
{\bf Representation dimension and Solomon zeta function}\\
\hspace{20mm} 
{\sc Osamu Iyama} 
\end{center}

Cline-Parshall-Scott introduced the concept of {\it quasi-hereditary algebras} (\S\ZBE) to study highest weight categories in the representation theory of Lie algebras and algebraic groups [CPS1,2]. Quasi-hereditary algebras were effectively applied in the representation theory of artin algebras as well by Dlab-Ringel [DR1,2,3] and many other authors. On the other hand, in the representation theory of orders, the concept of {\it overorders} and {\it overrings} (\S1.1), a non-commutative analogy of the normalization in the commutative ring theory, plays a crucial role. From an overring $\Gamma$ of an order $\Lambda$, we naturally obtain a full subcategory $\lat\Gamma$ of $\lat\Lambda$. Formulating this correspondence $\Gamma\mapsto\lat\Gamma$ categorically, we obtain the concept of the {\it rejection} (\S1,\S2). Recently it was effectively applied to study orders of finite representation type by the author [I1,2,3] and Rump [Ru1,2,3]. Originally Drozd-Kirichenko-Roiter found the one-point rejection (\S\ZAC) in their theory of Bass orders [DKR], and later Hijikata-Nishida applied the four-points rejection (\S\ZAE) to local orders of finite representation type and suggested a possibility of generalization [HN1,2,3].

In this paper, we will show that there exists a close relationship between quasi-hereditary algebras and the rejection from the viewpoint of the approximation theory of Auslander-Smalo [AS2]. As an application, we will solve two open problems [I4,5]. One concerns the {\it representation dimension} of artin algebras introduced by M. Auslander about 30 years ago [A1], and another concerns the {\it Solomon zeta functions} of orders introduced by L. Solomon about 25 years ago [S1,2]. It will turn out that the rejection relates these two quite different problems with each other closely.
\[\begin{array}{|c|c|}\hline
\mbox{orders (Krull dimension one)}&\mbox{artin algebras (Krull dimension zero)}\\ \hline
\mbox{overrings of an order $\Lambda$}&\mbox{factor algebras of an artin algebra $\Lambda$}\\
\parallel&\parallel\\
\mbox{rejective subcategories of $\lat\Lambda$}&\mbox{rejective subcategories of $\mod\Lambda$}\\
\cap&\cap\\
\mbox{right rejective subcategories of $\lat\Lambda$}&\mbox{right rejective subcategories of $\mod\Lambda$}\\ 
\downarrow&\downarrow\\
\mbox{Solomon's second conjecture}&\mbox{finiteness of representation dimension}\\ \hline
\end{array}\]

The representation theory of orders looks like that of artin algebras, and the key diagram above shows their correspondence. We will use the upper part in \S1, the middle part in \S2, the right part in \S3, and the left part in \S4. In \S5, we will relate Solomon zeta functions with Ringel-Hall algebras. Moreover, we will construct the universal enveloping algebra $U(\dn{gl}_n)$ and a family of irreducible $U(\dn{sl}_n)$-modules by using the Hall algebras of hereditary orders [I7]. This is closely related to Ringel's construction of the quantum group of type $\widetilde{A}_{n-1}$ [R2][Lu][Sc][As].

In the rest of this paper, any module is assumed to be a left module. For a ring $\Lambda$, we denote by $J_\Lambda$ the Jacobson radical of $\Lambda$, and by $\mod\Lambda$ (resp. $\pr\Lambda$) the category of finitely generated (resp. finitely generated projective) $\Lambda$-modules. For an artin algebra $\Lambda$, we denote by $(\ )^*:\mod\Lambda\leftrightarrow\mod\Lambda^{op}$ the duality. For an additive category $\cc$, we denote by $\ind\cc$ the set of isoclasses of indecomposable objects in $\cc$.

\vskip1em{\bf\ZA\ Overorders and the rejection}

Let $R$ be a commutative noetherian domain with the quotient field of $K$. An $R$-algebra $\Lambda$ is called an {\it $R$-order} if it is a finitely generated projective $R$-module. In this case, $\Lambda$ forms a subring of a finite dimensional $K$-algebra $A:=\Lambda\otimes_RK$. For an $R$-order $\Lambda$, a left $\Lambda$-module $L$ is called a {\it $\Lambda$-lattice} if it is a finitely generated projective $R$-module. We denote by $\lat\Lambda$ the category of $\Lambda$-lattices. Then $(\ )^*=\hom_R(\ ,R)$ gives a duality between $\lat\Lambda$ and $\lat\Lambda^{op}$. We call $\rin\Lambda:=(\pr\Lambda^{op})^*$ the category of {\it injective} $\Lambda$-lattices.

In the rest of this section, assume that $R$ is a complete discrete valuation ring. For example, $(R,K)$ is $(\zzz_p,\qqq_p)$ or $(k[[x]],k((x)))$ for a field $k$. Then $\lat\Lambda$ forms a {\it Krull-Schmidt category}, namely any object is isomorphic to a finite direct sum of objects whose endomorphism rings are local [CR]. Moreover, we assume that $A=\Lambda\otimes_RK$ is semisimple. Then the existence theorem of Auslander-Reiten sequences holds in $\lat\Lambda$ [A2][Ro2]. We denote by $\dn{A}(\Lambda)$ the Auslander-Reiten quiver of $\Lambda$ [ARS][Y], which is a quiver with the set of vertices $\ind(\lat\Lambda)$.

\vskip1em{\bf\ZAA\ Definition }
Let $\Gamma$ be another $R$-order. We call $\Gamma$ an {\it overorder} of $\Lambda$ if $A\supset\Gamma\supseteq\Lambda$ holds. More generally, we call $\Gamma$ an {\it overring} of $\Lambda$ if $A/I\supset\Gamma\supseteq(\Lambda+I)/I$ holds for some ideal $I$ of $A$. Then the natural morphism $\Lambda\rightarrow\Gamma$ induces a full faithful functor $\lat\Gamma\rightarrow\lat\Lambda$. Thus $\lat\Gamma$ can be regarded as a full subcategory of $\lat\Lambda$, and $\ind(\lat\Gamma)$ forms a subset of $\ind(\lat\Lambda)$. We denote by $\ovo{\Lambda}$ (resp. $\ovr{\Lambda}$) the set of overorders (resp. overrings) of $\Lambda$. We introduce a partial order $\subseteq$ on $\ovr{\Lambda}$ as follows: For $\Gamma_i\in\ovr{\Lambda}$ ($i=1,2$), define $\Gamma_1\subseteq\Gamma_2$ if and only if $\Gamma_1^\prime\subseteq\Gamma_2^\prime$ holds as subsets of $A$ for the pull-back $\Gamma_i^\prime$ of $\Gamma_i$ to $A$. Then the correspondence $\Gamma\mapsto\ind(\lat\Gamma)$ gives an injection $\ovr{\Lambda}\to 2^{\ind(\lat\Lambda)}$, which reverses the partial order $\subseteq$ [I3]. An $R$-order $\Lambda$ with $\ovo{\Lambda}=\{\Lambda\}$ is called {\it maximal}.

\vskip1em{\bf\ZAB\ Example }
(1) Put $\Lambda:=\left(\begin{array}{cc}{\scriptstyle R}&{\scriptstyle R}\\ {\scriptstyle J_R^n}&{\scriptstyle R}\end{array}\right)\subset A:=\ma_2(K)$ ($n\ge0$). Then $\ovr{\Lambda}=\ovo{\Lambda}$ is given by the following. (For simplicity, put $n=2$.)
{\scriptsize\[
\begin{array}{ccccc}
\left(\begin{array}{cc}
{\scriptstyle R}&{\scriptstyle J_R^{-2}}\\ {\scriptstyle J_R^2}&{\scriptstyle R}\end{array}\right)&&&&\\
\cup&&&&\\
\left(\begin{array}{cc}{\scriptstyle R}&{\scriptstyle J_R^{-1}}\\ {\scriptstyle J_R^2}&{\scriptstyle R}\end{array}\right)&\subset&\left(\begin{array}{cc}{\scriptstyle R}&{\scriptstyle J_R^{-1}}\\ {\scriptstyle J_R}&{\scriptstyle R}\end{array}\right)&&\\
\cup&&\cup&&\\
\Lambda=\left(\begin{array}{cc}{\scriptstyle R}&{\scriptstyle R}\\ {\scriptstyle J_R^2}&{\scriptstyle R}\end{array}\right)&\subset&\left(\begin{array}{cc}{\scriptstyle R}&{\scriptstyle R}\\ {\scriptstyle J_R}&{\scriptstyle R}\end{array}\right)&\subset&\left(\begin{array}{cc}{\scriptstyle R}&{\scriptstyle R}\\ {\scriptstyle R}&{\scriptstyle R}\end{array}\right)
\end{array}\]}

The corresponding subsets of $\ind(\lat\Lambda)$ are given by the following.
{\scriptsize\[\begin{array}{ccccc}
\left\{{R\choose J_R^2}\right\}&&&&\\
\cap&&&&\\
\left\{{R\choose J_R^2},{R\choose J_R}\right\}&\supset&\left\{{R\choose J_R}\right\}&&\\
\cap&&\cap&&\\
\left\{{R\choose J_R^2},{R\choose J_R},{R\choose R}\right\}&\supset&\left\{{R\choose J_R},{R\choose R}\right\}&\supset&\left\{{R\choose R}\right\}
\end{array}\]}

(2) Put $\Lambda=\Lambda_n:=\{(x,y)\in R\times R\ |\ x-y\in J_R^n\}\subset A:=K\times K$ ($n\ge0$). Then $\ovr{\Lambda}=\ovo{\Lambda}\cup\{ 0\times R,R\times0\}$ is given by the following.
\[\begin{array}{ccc}
&0\times R&\\
&\cup&\\
\Lambda_n\subset\Lambda_{n-1}\subset\cdots\subset\Lambda_2\subset\Lambda_1\subset&\Lambda_0=R\times R&\subset R\times 0
\end{array}\]

The corresponding subsets of $\ind(\lat\Lambda)$ are given by the following.
{\scriptsize\[\begin{array}{ccc}
&\{ 0\times R\}&\\
&\cap&\\
\left\{\begin{array}{c}
{ \Lambda_i}\\ {\scriptstyle (1\le i\le n)}\\ { R\times0}\\ { 0\times R}\end{array}\right\}\supset
\left\{\begin{array}{c}{ \Lambda_i}\\ {\scriptstyle (1\le i\le n-1)}\\ { R\times 0}\\ { 0\times R}\end{array}\right\}\supset\cdots\supset
\left\{\begin{array}{c}{ \Lambda_1,\Lambda_2}\\ { R\times0}\\ { 0\times R}\end{array}\right\}\supset
\left\{\begin{array}{c}{ \Lambda_1}\\ { R\times0}\\ { 0\times R}\end{array}\right\}\supset
&\left\{\begin{array}{c}{ R\times0}\\ { 0\times R}\end{array}\right\}&\\
&\cup&\\
&\{ R\times0\}&
\end{array}\]}

We describe their Auslander-Reiten quivers, which look like the Dynkin diagram $\ddd_i$.
{\scriptsize\begin{eqnarray*}&\begin{array}{ccr}
\dn{A}(\Lambda_n)&&\begin{array}{ccc}
&\stackrel{0\times R}{\bullet}&\\
&\downarrow\uparrow&\\
\stackrel{\Lambda_n}{\bullet}
\def\arraystretch{.2}\begin{array}{cc}\longrightarrow\\
\longleftarrow\end{array}
\stackrel{\Lambda_{n-1}}{\bullet}
\def\arraystretch{.2}\begin{array}{cc}\longrightarrow\\
\longleftarrow\end{array}\cdots
\def\arraystretch{.2}\begin{array}{cc}\longrightarrow\\
\longleftarrow\end{array}
\stackrel{\Lambda_2}{\bullet}
\def\arraystretch{.2}\begin{array}{cc}\longrightarrow\\
\longleftarrow\end{array}
&\stackrel{\Lambda_1}{\bullet}&
\def\arraystretch{.2}\begin{array}{cc}\longrightarrow\\
\longleftarrow\end{array}
\stackrel{R\times 0}{\bullet}
\end{array}\\
\dn{A}(\Lambda_{n-1})&&\begin{array}{ccc}
&\stackrel{0\times R}{\bullet}&\\
&\downarrow\uparrow&\\
\stackrel{\Lambda_{n-1}}{\bullet}
\def\arraystretch{.2}\begin{array}{cc}\longrightarrow\\
\longleftarrow\end{array}\cdots
\def\arraystretch{.2}\begin{array}{cc}\longrightarrow\\
\longleftarrow\end{array}
\stackrel{\Lambda_2}{\bullet}
\def\arraystretch{.2}\begin{array}{cc}\longrightarrow\\
\longleftarrow\end{array}
&\stackrel{\Lambda_1}{\bullet}&
\def\arraystretch{.2}\begin{array}{cc}\longrightarrow\\
\longleftarrow\end{array}
\stackrel{R\times 0}{\bullet}
\end{array}\\
\cdots&&\cdots\\
\dn{A}(\Lambda_2)&&\begin{array}{ccc}
&\stackrel{0\times R}{\bullet}&\\
&\downarrow\uparrow&\\
\stackrel{\Lambda_2}{\bullet}
\def\arraystretch{.2}\begin{array}{cc}\longrightarrow\\
\longleftarrow\end{array}
&\stackrel{\Lambda_1}{\bullet}&
\def\arraystretch{.2}\begin{array}{cc}\longrightarrow\\
\longleftarrow\end{array}
\stackrel{R\times 0}{\bullet}
\end{array}\\
\dn{A}(\Lambda_{1})&&\begin{array}{ccc}
&\stackrel{0\times R}{\bullet}&\\
&\downarrow\uparrow&\\
&\stackrel{\Lambda_1}{\bullet}&
\def\arraystretch{.2}\begin{array}{cc}\longrightarrow\\
\longleftarrow\end{array}
\stackrel{R\times 0}{\bullet}
\end{array}\\
\end{array}&\\
&\dn{A}(\Lambda_{0})=\left(\begin{array}{cc}
\stackrel{0\times R}{\bullet}&
\stackrel{R\times 0}{\bullet}\\
\carrow&\carrow
\end{array}\right),\ \ \ \ \ 
\dn{A}(R\times 0)=\left(\begin{array}{c}
\stackrel{R\times 0}{\bullet}\\
\carrow
\end{array}\right),\ \ \ \ \ 
\dn{A}(0\times R)=\left(\begin{array}{c}
\stackrel{0\times R}{\bullet}\\
\carrow
\end{array}\right)&
\end{eqnarray*}}

(3) Put $R:=k[[x^2]]$ and $\Lambda=\Lambda_n:=R+Rx^{2n+1}\subset A:=k((x))$ ($n\ge0$). Then $\ovr{\Lambda}=\ovo{\Lambda}$ is given by the following.
\[\Lambda_n\subset\Lambda_{n-1}\subset\cdots\subset\Lambda_2\subset\Lambda_1\subset\Lambda_0\]

The corresponding subsets of $\ind(\lat\Lambda)$ are given by the following.
\[\{\Lambda_i\}_{0\le i\le n}\supset\{\Lambda_i\}_{0\le i\le n-1}\supset\cdots\supset\{\Lambda_0,\Lambda_1,\Lambda_2\}\supset\{\Lambda_0,\Lambda_1\}\supset\{\Lambda_0\}\]

We describe their Auslander-Reiten quivers, which look like the Dynkin diagram $\aaa_i$.
{\scriptsize\[\begin{array}{ccr}
\dn{A}(\Lambda_n)&&
\stackrel{\Lambda_n}{\bullet}
\def\arraystretch{.2}\begin{array}{cc}\longrightarrow\\
\longleftarrow\end{array}
\stackrel{\Lambda_{n-1}}{\bullet}
\def\arraystretch{.2}\begin{array}{cc}\longrightarrow\\
\longleftarrow\end{array}\cdots
\def\arraystretch{.2}\begin{array}{cc}\longrightarrow\\
\longleftarrow\end{array}
\stackrel{\Lambda_2}{\bullet}
\def\arraystretch{.2}\begin{array}{cc}\longrightarrow\\
\longleftarrow\end{array}
\stackrel{\Lambda_1}{\bullet}
\def\arraystretch{.2}\begin{array}{cc}\longrightarrow\\
\longleftarrow\end{array}
\stackrel{\Lambda_0}{\bullet}\carrow\\
\dn{A}(\Lambda_{n-1})&&
\stackrel{\Lambda_{n-1}}{\bullet}
\def\arraystretch{.2}\begin{array}{cc}\longrightarrow\\
\longleftarrow\end{array}\cdots
\def\arraystretch{.2}\begin{array}{cc}\longrightarrow\\
\longleftarrow\end{array}
\stackrel{\Lambda_2}{\bullet}
\def\arraystretch{.2}\begin{array}{cc}\longrightarrow\\
\longleftarrow\end{array}
\stackrel{\Lambda_1}{\bullet}
\def\arraystretch{.2}\begin{array}{cc}\longrightarrow\\
\longleftarrow\end{array}
\stackrel{\Lambda_0}{\bullet}\carrow\\
\cdots&&\cdots\\
\dn{A}(\Lambda_{1})&&
\stackrel{\Lambda_1}{\bullet}
\def\arraystretch{.2}\begin{array}{cc}\longrightarrow\\
\longleftarrow\end{array}
\stackrel{\Lambda_0}{\bullet}\carrow\\
\dn{A}(\Lambda_{0})&&
\stackrel{\Lambda_0}{\bullet}\carrow
\end{array}\]}

In (2) and (3) above, it is remarkable that one-point is removed repeatedly from $\dn{A}(\Lambda_n)$ to $\dn{A}(\Lambda_0)$. Similarly in (1) above, any consecutive two orders $\Gamma_2\subset\Gamma_1$ satisfies $\#(\ind(\lat\Gamma_2)-\ind(\lat\Gamma_1))=1$. Although such a phenomenon does not hold in general, we will analyze it by introducing the following concept.

\vskip1em{\bf\ZAC\ }
We call a subset $S$ of $\ind(\lat\Lambda)$ {\it rejectable} (resp. {\it strictly rejectable}) if there exists an overring (resp. overorder) $\Gamma$ of $\Lambda$ such that $S=\ind(\lat\Lambda)-\ind(\lat\Gamma)$ [I1,2,3][Ru1,2,3]. Then the correspondence $\Gamma\mapsto\ind(\lat\Lambda)-\ind(\lat\Gamma)$ gives a bijection from $\ovr{\Lambda}$ to the set of rejectable subsets of $\ind(\lat\Lambda)$, which preserves the partial order $\subseteq$. We can explain the phenomenon above by the following lemma of Drozd-Kirichenko [DK1][HN2].

\vskip1em{\bf (One-point rejection) }{\it
Let $\Lambda$ be an order and $X\in\ind(\lat\Lambda)$. Then a singleton set $\{ X\}$ is rejectable if and only if $X\in\pr\Lambda\cap\rin\Lambda$.}

\vskip1em
This is fundamental in the theory of Bass orders [DKR][Ro1][HN1,2]. Recall that an order $\Lambda$ is called {\it Gorenstein} if $\Lambda\in\rin\Lambda$, and {\it Bass} if any overorder of $\Lambda$ is Gorenstein. Three examples in \ZAB\ are Bass orders. For any Bass order $\Lambda$, the one-point rejection assures the existence of a sequence of overorders $\Lambda=\Lambda_n\subset\Lambda_{n-1}\subset\cdots\subset\Lambda_0$ with $\#(\ind(\lat\Lambda_i)-\ind(\lat\Lambda_{i-1}))=1$ for any $i$. In particular, it is shown that $\dn{A}(\Lambda)$ is `given by' a Dynkin diagram (e.g. [W]). Here we won't go further on Bass orders.

\vskip1em{\bf\ZAD\ Example }
Let us observe rejectable subsets for non-Bass orders. Let $\Lambda$ be the completion of a simple curve singularity [Y]. Then $\Lambda$ forms an order over a complete discrete valuation ring $R$ [DW]. For the case $\aaa_n:k[[X,Y]]/(X^{n+1}+Y^2)$, $\Lambda$ is shown to be a Bass order given in \ZAB(2) (if $n$ is odd) or \ZAB(3) (if $n$ is even). Thus we will study the case $\ddd_n:k[[X,Y]]/(X^{n-1}+XY^2)$, and we assume that $n$ is an odd $2m+1$ (the even case is quite similar). Then $\ovr{\Lambda}=\ovo{\Lambda}\cup\{\Omega,\Gamma_i\}_{1\le i\le m}$ is the diagram below, where $\Omega$ is a maximal order, $\Gamma_i$ is a simple curve singularity of type $\aaa_{2i-2}$, and $\Lambda_i$ is a unique minimal overorder of a simple curve singularity of type $\ddd_{2i+1}$.
{\scriptsize\[\begin{array}{ccccccccccccc}
&&\Gamma_m&\subset&\Gamma_{m-1}&\subset&\cdots&\subset&\Gamma_2&\subset&\Gamma_1&&\\
&&\cup&&\cup&&&&\cup&&\cup&&\\
&&\Omega\times\Gamma_m&\subset&\Omega\times\Gamma_{m-1}&\subset&\cdots&\subset&\Omega\times\Gamma_2&\subset&\Omega\times\Gamma_1&\subset&\Omega\\
&&\cup&&\cup&&&&\cup&&\cup&&\\
\Lambda&\subset&\Lambda_m&\subset&\Lambda_{m-1}&\subset&\cdots&\subset&\Lambda_2&\subset&\Lambda_1&&
\end{array}\]}

Then $\dn{A}(\Lambda)$ is the quiver below for some $X_i\in\lat\Lambda$, where the bottom vertices are identified with the top vertices and $\Lambda_1^*=\Lambda_1$ holds.
{\scriptsize
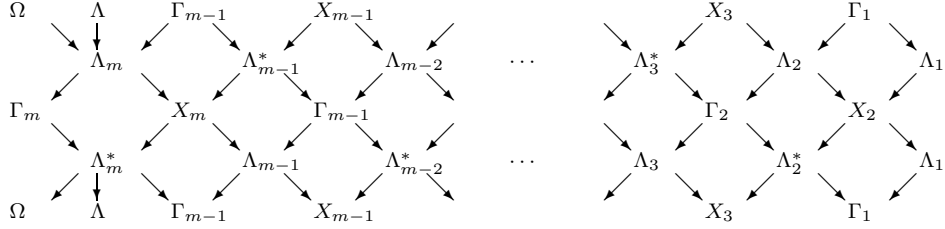
\[
\begin{array}{cccccccccccccccccccccccccccc}
\Q{\Omega}&&\Q{\Lambda}&&\Q{\Gamma_{m-1}}&&&&\Q{X_{m-1}}&&&&\Q{}&&&&\Q{}&&&&\Q{X_3}&&&&\Q{\Gamma_1}&&\\
&\QDR&\QD&\QDL&&\QDR&&\QDL&&\QDR&&\QDL&&&&&&\QDR&&\QDL&&\QDR&&\QDL&&\QDR&\\
&&\Q{\Lambda_m}&&&&\Q{\Lambda_{m-1}^*}&&&&\Q{\Lambda_{m-2}}&&&&\Q{\cdots}&&&&\Q{\Lambda_3^*}&&&&\Q{\Lambda_2}&&&&\Q{\Lambda_1}\\
&\QDL&&\QDR&&\QDL&&\QDR&&\QDL&&\QDR&&&&&&\QDL&&\QDR&&\QDL&&\QDR&&\QDL&\\
\Q{\Gamma_m}&&&&\Q{X_m}&&&&\Q{\Gamma_{m-1}}&&&&\Q{}&&&&\Q{}&&&&\Q{\Gamma_2}&&&&\Q{X_2}&&\\
&\QDR&&\QDL&&\QDR&&\QDL&&\QDR&&\QDL&&&&&&\QDR&&\QDL&&\QDR&&\QDL&&\QDR&\\
&&\Q{\Lambda_m^*}&&&&\Q{\Lambda_{m-1}}&&&&\Q{\Lambda_{m-2}^*}&&&&\Q{\cdots}&&&&\Q{\Lambda_3}&&&&\Q{\Lambda_2^*}&&&&\Q{\Lambda_1}\\
&\QDL&\QD&\QDR&&\QDL&&\QDR&&\QDL&&\QDR&&&&&&\QDL&&\QDR&&\QDL&&\QDR&&\QDL&\\
\Q{\Omega}&&\Q{\Lambda}&&\Q{\Gamma_{m-1}}&&&&\Q{X_{m-1}}&&&&\Q{}&&&&\Q{}&&&&\Q{X_3}&&&&\Q{\Gamma_1}&&
\end{array}\]}

Let us describe $\dn{A}(\Lambda_i)$ for any $i$. Since we can apply the one-point rejection \ZAC\ to $\dn{A}(\Lambda)$, we obtain $\dn{A}(\Lambda_m)$ below by removing $\{\Lambda\}$ from $\dn{A}(\Lambda)$.
{\scriptsize
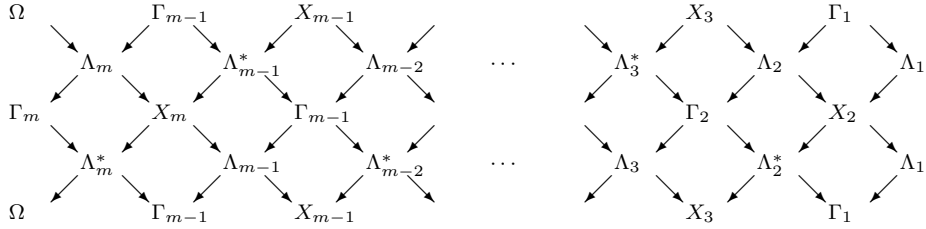
\[
\begin{array}{cccccccccccccccccccccccccccc}
\Q{\Omega}&&&&\Q{\Gamma_{m-1}}&&&&\Q{X_{m-1}}&&&&\Q{}&&&&\Q{}&&&&\Q{X_3}&&&&\Q{\Gamma_1}&&\\
&\QDR&&\QDL&&\QDR&&\QDL&&\QDR&&\QDL&&&&&&\QDR&&\QDL&&\QDR&&\QDL&&\QDR&\\
&&\Q{\Lambda_m}&&&&\Q{\Lambda_{m-1}^*}&&&&\Q{\Lambda_{m-2}}&&&&\Q{\cdots}&&&&\Q{\Lambda_3^*}&&&&\Q{\Lambda_2}&&&&\Q{\Lambda_1}\\
&\QDL&&\QDR&&\QDL&&\QDR&&\QDL&&\QDR&&&&&&\QDL&&\QDR&&\QDL&&\QDR&&\QDL&\\
\Q{\Gamma_m}&&&&\Q{X_m}&&&&\Q{\Gamma_{m-1}}&&&&\Q{}&&&&\Q{}&&&&\Q{\Gamma_2}&&&&\Q{X_2}&&\\
&\QDR&&\QDL&&\QDR&&\QDL&&\QDR&&\QDL&&&&&&\QDR&&\QDL&&\QDR&&\QDL&&\QDR&\\
&&\Q{\Lambda_m^*}&&&&\Q{\Lambda_{m-1}}&&&&\Q{\Lambda_{m-2}^*}&&&&\Q{\cdots}&&&&\Q{\Lambda_3}&&&&\Q{\Lambda_2^*}&&&&\Q{\Lambda_1}\\
&\QDL&&\QDR&&\QDL&&\QDR&&\QDL&&\QDR&&&&&&\QDL&&\QDR&&\QDL&&\QDR&&\QDL&\\
\Q{\Omega}&&&&\Q{\Gamma_{m-1}}&&&&\Q{X_{m-1}}&&&&\Q{}&&&&\Q{}&&&&\Q{X_3}&&&&\Q{\Gamma_1}&&
\end{array}\]}

We obtain $\dn{A}(\Lambda_{m-1})$ below by removing four-points
{\scriptsize$\begin{array}{ccc}
\Lambda_m&\rightarrow&X_m\\
\downarrow&&\downarrow\\
\Gamma_m&\rightarrow&\Lambda_m^*
\end{array}$} from $\dn{A}(\Lambda_m)$.
{\scriptsize\[\dn{A}(\Lambda_{m-1})\ \ \ \begin{array}{cccccccccccccccccccccccc}
\Q{\Gamma_{m-1}}&&&&\Q{X_{m-1}}&&&&\Q{}&&&&\Q{}&&&&\Q{X_3}&&&&\Q{\Gamma_1}&&\\
&\QDR&&\QDL&&\QDR&&\QDL&&&&&&\QDR&&\QDL&&\QDR&&\QDL&&\QDR&\\
&&\Q{\Lambda_{m-1}^*}&&&&\Q{\Lambda_{m-2}}&&&&\Q{\cdots}&&&&\Q{\Lambda_3^*}&&&&\Q{\Lambda_2}&&&&\Q{\Lambda_1}\\
&\QDL&&\QDR&&\QDL&&\QDR&&&&&&\QDL&&\QDR&&\QDL&&\QDR&&\QDL&\\
\Q{\Omega}&&&&\Q{\Gamma_{m-1}}&&&&\Q{}&&&&\Q{}&&&&\Q{\Gamma_2}&&&&\Q{X_2}&&\\
&\QDR&&\QDL&&\QDR&&\QDL&&&&&&\QDR&&\QDL&&\QDR&&\QDL&&\QDR&\\
&&\Q{\Lambda_{m-1}}&&&&\Q{\Lambda_{m-2}^*}&&&&\Q{\cdots}&&&&\Q{\Lambda_3}&&&&\Q{\Lambda_2^*}&&&&\Q{\Lambda_1}\\
&\QDL&&\QDR&&\QDL&&\QDR&&&&&&\QDL&&\QDR&&\QDL&&\QDR&&\QDL&\\
\Q{\Gamma_{m-1}}&&&&\Q{X_{m-1}}&&&&\Q{}&&&&\Q{}&&&&\Q{X_3}&&&&\Q{\Gamma_1}&&
\end{array}\]}

Similarly, we obtain $\dn{A}(\Lambda_{m-2})$ below by removing four-points
{\scriptsize$\begin{array}{ccc}
\Lambda_{m-1}&\rightarrow&X_{m-1}\\
\downarrow&&\downarrow\\
\Gamma_{m-1}&\rightarrow&\Lambda_{m-1}^*
\end{array}$} from $\dn{A}(\Lambda_{m-1})$. Thus four points
{\scriptsize$\begin{array}{ccc}
\Lambda_i&\rightarrow&X_i\\
\downarrow&&\downarrow\\
\Gamma_i&\rightarrow&\Lambda_i^*
\end{array}$}
are removed repeatedly, and finally we obtain $\dn{A}(\Lambda_{2})$ and $\dn{A}(\Lambda_{1})$ below.
{\scriptsize\[\dn{A}(\Lambda_{2})\ \ \ \begin{array}{ccccccc}
\Q{\Omega}&&&&\Q{\Gamma_1}&&\\
&\QDR&&\QDL&&\QDR&\\
&&\Q{\Lambda_2}&&&&\Q{\Lambda_1}\\
&\QDL&&\QDR&&\QDL&\\
\Q{\Gamma_2}&&&&\Q{X_2}&&\\
&\QDR&&\QDL&&\QDR&\\
&&\Q{\Lambda_2^*}&&&&\Q{\Lambda_1}\\
&\QDL&&\QDR&&\QDL&\\
\Q{\Omega}&&&&\Q{\Gamma_1}&&
\end{array}\ \ \ \ \ \ \ \ \ \ \ \ \ \ \ \ \ \ \ \ \ \ \ \ \ 
\dn{A}(\Lambda_{1})\ \ \ \begin{array}{ccccccc}
\Q{\Gamma_1}&&\\
&\QDR&\\
&&\Q{\Lambda_1}\\
&\QDL&\\
\Q{\Omega}&&\\
&\QDR&\\
&&\Q{\Lambda_1}\\
&\QDL&\\
\Q{\Gamma_1}&&
\end{array}\]}

\vskip0em{\bf\ZAE\ }
It is natural to ask why four-points are removed repeatedly in \ZAD. The following analogy of the one-point rejection in \ZAC\ gives a reason.

\vskip1em{\bf (Four-points rejection) }{\it Let $\Lambda$ be an order, $\tau$ the Auslander translate and $S=\{ P,X,Y,I\}$ a subset of $\ind(\lat\Lambda)$. If $S$ has the form {\scriptsize$\begin{array}{ccc}
P&\rightarrow&X\\
\downarrow&&\downarrow\\
Y&\rightarrow&I
\end{array}$} in $\dn{A}(\Lambda)$ with $P\in\pr\Lambda$, $I\in\rin\Lambda$ and $P=\tau I$, then $S$ is rejectable.}

\vskip1em

As the one-point rejection is fundamental for Bass orders, the four-points rejection is fundamental for local orders of fnite representation type [DK2][HN3]. More generally, the {\it general rejection} for any finite subset $S$ of $\ind(\lat\Lambda)$, a combinatorial condition which decides whether $S$ is rejectable or not, is given by the author [I1,2,3]. It depends only on the restriction of the quiver $\dn{A}(\Lambda)$ to $S$, so the rejectability of $S$ is independent of the outside of $S$. Although we omit it here, we will give examples in \ZAF\ below.

\vskip1em{\bf\ZAF\ Example }[I1]
Let $\Lambda$ be an order, $\tau$ the Auslander translate and $S$ a subset of $\ind(\lat\Lambda)$. We call $S$ {\it minimal rejectable} if all rejectable subsets of $\ind(\lat\Lambda)$ contained in $S$ are $S$ and $\emptyset$. Then one can show that $S\cap\pr\Lambda$ and $S\cap\rin\Lambda$ are singleton sets.

Now we assume that $S$ is a subset of $\ind(\lat\Lambda)$ with $\# S\le 4$. Then we can show that $S$ is minimal strictly rejectable if and only if $S$ has one of the forms (1)--(12) below, where we put $\{ P\}:=S\cap\pr\Lambda$ and $\{ I\}:=S\cap\rin\Lambda$. The case (1) is the one-point rejection \ZAC, and the case (12) is the four-points rejection \ZAE.
{\scriptsize\begin{eqnarray*}
(1)&\stackrel{P=I}{\bullet}& \\
(2)&\stackrel{P}{\bullet}\longrightarrow\stackrel{I}{\bullet}&\\
(3)&\stackrel{P}{\bullet}\stackrel{(a\ b)}{\longrightarrow}\bullet\longrightarrow\stackrel{I}{\bullet}&ab\le2\\
(4)&\stackrel{P}{\bullet}\longrightarrow\bullet\stackrel{(a\ b)}{\longrightarrow}\stackrel{I}{\bullet}&ab\le2\\
(5)&\stackrel{P}{\bullet}\stackrel{(a\ b)}{\longrightarrow}\bullet\stackrel{(b\ a)}{\longrightarrow}\stackrel{I}{\bullet}&P=\tau I, ab\le3\\
(6)&\stackrel{P}{\bullet}\stackrel{(a\ b)}{\longrightarrow}\bullet\longrightarrow\bullet\longrightarrow\stackrel{I}{\bullet}&ab\le2\\
(7)&\stackrel{P}{\bullet}\longrightarrow\bullet\stackrel{(a\ b)}{\longrightarrow}\bullet\longrightarrow\stackrel{I}{\bullet}&ab\le2\\
(8)&\stackrel{P}{\bullet}\longrightarrow\bullet\longrightarrow\bullet\stackrel{(a\ b)}{\longrightarrow}\stackrel{I}{\bullet}&ab\le2\\
(9)&\stackrel{P}{\bullet}\stackrel{(a\ b)}{\longrightarrow}\bullet\stackrel{(b\ a)}{\longrightarrow}\stackrel{X}{\bullet}\longrightarrow\stackrel{I}{\bullet}&P=\tau X, ab\le3\\
(10)&\stackrel{P}{\bullet}\longrightarrow\stackrel{X}{\bullet}\stackrel{(a\ b)}{\longrightarrow}\bullet\stackrel{(b\ a)}{\longrightarrow}\stackrel{I}{\bullet}&X=\tau I, ab\le3\\
(11)&\begin{array}{ccc}
\stackrel{P}{\bullet}\longrightarrow&\bullet&\longrightarrow\stackrel{I}{\bullet}\\
&\downarrow\uparrow&\\
&{\scriptstyle X}\bullet&
\end{array}&P=\tau I,X=\tau X\\
(12)&\begin{array}{ccc}
\ \ \bullet&\longrightarrow&\stackrel{I}{\bullet}\\
\ \ \uparrow&&\uparrow\\
{\scriptstyle P}\bullet&\longrightarrow&\bullet
\end{array}&P=\tau I
\end{eqnarray*}}

\vskip0em{\bf\ZAG\ }
Let us consider artin algebras. Let $\Lambda$ be an artin algebra and $S$ a subset of $\ind(\mod\Lambda)$. We call $S$ {\it rejectable} if there exists a factor algebra $\Gamma=\Lambda/I$ of $\Lambda$ such that $S=\ind(\mod\Lambda)-\ind(\mod\Gamma)$. Then the general rejection in \ZAE\ (including the one-point and four-points rejections) works for artin algebras as well. This observation motivates the categorical formulation of the rejection in the following section \S2.

\vskip1em
{\bf\ZB\ Approximation, rejective subcategories and quasi-hereditary algebras }

Let $\cc$ be an additive category, $\cc(X,Y):=\hom_{\cc}(X,Y)$, and $fg$ the composition of $f\in\cc(X,Y)$ and $g\in\cc(Y,Z)$. Throughout this paper, {\it any subcategory is assumed to be full and closed under isomorphisms, direct sums and direct summands.} For any $X\in\cc$, we denote by $\add X$ the full subactegory of $\cc$ consisting of direct summands of a finite direct sum of $X$. We call $X\in\cc$ an {\it additive generator} of $\cc$ if $\add X=\cc$ holds. In this section, we denote by $R$ a complete discrete valuation ring.

\vskip1em{\bf\ZBA\ } Let $I$ be an {\it ideal} of $\cc$, namely a subgroup $I(X,Y)$ of $\cc(X,Y)$ is given for any $X,Y\in\cc$ such that $fgh\in I(W,Z)$ holds for any $f\in\cc(W,X)$, $g\in I(X,Y)$ and $h\in\cc(Y,Z)$. Then the {\it factor category} $\cc/I$ is defined by $\ob(\cc)=\ob(\cc/I)$ and $(\cc/I)(X,Y):=\cc(X,Y)/I(X,Y)$ for any $X,Y\in\cc$.

We call $f\in I(Y,X)$ a {\it right $I$-approximation} of $X$ if $\cc(\ ,Y)\stackrel{\cdot f}{\rightarrow}I(\ ,X)\rightarrow0$ is exact. We call $I$ {\it right finitely generated} if any $X\in\cc$ has a right $I$-approximation. Dually, a {\it left $I$-approximation} and a {\it left finitely generated} ideal are defined.

The concept of the approximation is introduced by Auslander-Smalo [AS2], and there are many important examples. We denote by $J_{\cc}$ the {\it Jacobson radical} of $\cc$, namely $J_{\cc}$ is the ideal of $\cc$ such that $J_{\cc}(X,X)$ forms the usual Jacobson radical of the ring $\cc(X,X)$ for any $X\in\cc$. For a subcategory $\cc^\prime$ of $\cc$, we denote by $[\cc^\prime]$ the ideal of $\cc$ consisting of morphisms which factor through some object of $\cc^\prime$.

\vskip1em{\bf\ZBB\ Example }
(1) Let $\Lambda$ be a commutative local noetherian ring with a dualizing module, and $\cm\Lambda$ the category of maximal Cohen-Macaulay $\Lambda$-modules. Then the ideal $[\cm\Lambda]$ of $\mod\Lambda$ is left and right finitely generated by Auslander-Buchweitz theory [AB].

(2) For any $R$-order (resp. artin algebra) $\Lambda$, the ideal $J_{\lat\Lambda}$ of $\lat\Lambda$ (resp. $J_{\mod\Lambda}$ of $\mod\Lambda$) is left and right finitely generated by Auslander-Reiten theory [ARS][A2][Ro2].

\vskip1em{\bf\ZBC\ Definition }
Let us introduce a special class of right finitely generated subcategories. We call a subcategory $\cc^\prime$ of an additive category $\cc$ {\it right rejective} if any $X\in\cc$ has a right $[\cc^\prime]$-approximation $f\in\cc(Y,X)$ which is a monomorphism. In this case, $Y\in\cc^\prime$ holds. Dually, a {\it left rejective} subcategory is defined, and we call a left and right rejective subcategory a {\it rejective} subcategory. These definitions are equivalent to those in [I3,5]. 

In the sense of \ZBD\ below, we can regard the concept of rejective subcategories as a categorical formulation of the concept of overrings. But right rejective subcategories do not have such a representation theoretic meaning. There are much more right rejective subcategories than rejective subcategories, and we will apply them in \S3 and \S4. Notice that, for a class of additive categories $\cc$ called {\it $\tau$-categories}, we can characterize rejective subcategories of $\cc$ in terms of the factor category $\cc/[\cc^\prime]$ [I3]. This immediately gives the general rejection stated in \ZAE.

\vskip.5em{\bf\ZBCA\ }([I5]2.1.1) We collect basic facts. Let $\cc^\prime$ be a right rejective subcategory of $\cc$ and $\cc^{\prime\prime}$ a subcategory of $\cc^\prime$.

(1) Then $\cc^\prime/[\cc^{\prime\prime}]$ is a right rejective subcategory of $\cc/[\cc^{\prime\prime}]$. 

(2) If any $X\in\cc^\prime$ has a right $[\cc^{\prime\prime}]$-approximation which is a monomorphism {\it in $\cc$}, then $\cc^{\prime\prime}$ is a right rejective subcategory of $\cc$.

\vskip1em{\bf\ZBD\ Proposition }{\it[I3]
(1) Let $\Lambda$ be an $R$-order and $\cc:=\lat\Lambda$. Then a subcategory $\cc^\prime$ of $\cc$ is rejective if and only if $\cc^\prime=\lat\Gamma$ for an overring $\Gamma$ of $\Lambda$.

(2) Let $\Lambda$ be an artin algebra and $\cc:=\mod\Lambda$. Then a subcategory $\cc^\prime$ of $\cc$ is rejective if and only if $\cc^\prime=\mod\Gamma$ for a factor algebra $\Gamma=\Lambda/I$ of $\Lambda$.

In both cases, a right $[\cc^\prime]$-approximation and a left $[\cc^\prime]$-approximation of $X\in\cc$ are given by $\hom_\Lambda(\Gamma,X)\rightarrow X$ and $X\rightarrow(\Gamma\otimes_\Lambda X)^{**}$ respectively.}

\vskip1em{\bf\ZBE\ }
We recall quasi-hereditary algebras of Cline-Parshall-Scott [CPS1].

Let $\Lambda$ be an artin algebra and $I$ a two-sided ideal of $\Lambda$.
We call $I$ a {\it heredity ideal} of $\Lambda$ if $I^2=I\in\pr\Lambda$ and $IJ_\Lambda I=0$ hold. We call $\Lambda$ a {\it quasi-hereditary algebra} if there exists a chain $0=I_m\subseteq I_{m-1}\subseteq\cdots\subseteq I_0=\Lambda$ of ideals such that $I_{n-1}/I_n$ is a heredity ideal of $\Lambda/I_n$ for any $n$ ($0<n\le m$). In this case, $\gl\Lambda\le 2m-2$ holds [DR1]. Moreover, $\Lambda$ is a quasi-hereditary algebra if and only if $\mod\Lambda$ forms a {\it highest weight category} [CPS1].

We recall an one-dimensional analogy of quasi-hereditary algebras introduced by K\"onig-Wiedemann [KW]. An $R$-order $\Lambda$ is called a {\it quasi-hereditary order} if there exists an idempotent $e$ of $\Lambda$ such that $e\Lambda e$ is a maximal order and $\Lambda/\Lambda e\Lambda$ is a quasi-hereditary algebra defined above. Then $\gl\Lambda<\infty$ holds again.

\vskip1em{\bf\ZBF\ Definition }
Let $\cc$ be an additive category. A chain $\cc^\prime=\cc_m\subseteq\cc_{m-1}\subseteq\cdots\subseteq\cc_0=\cc$ of subcategories is called a {\it right rejective chain from $\cc$ to $\cc^\prime$} if $J_{\cc_n/[\cc_{n+1}]}=0$ holds and $\cc_{n+1}$ is a right rejective subcategory of $\cc$ for any $n$ ($0\le n<m$) [I5]\footnote{In this case, $\cc_{n+1}$ is a right rejective subcategory of $\cc_n$ for any $n$. Thus the definition above is slightly stronger than that in [I5]. Since the definition above was essentially used in [I5] as well, we adopt it. (We note here that the latter assertion in [I5]2.1.1 should be the following: {\it If $\cc^{\prime\prime}$ is a right rejective subcategory of $\cc^\prime$ with a counit $\epsilon^\prime$ such that $\epsilon^\prime_X$ is a monomorphism in $\cc$ for any $X\in\cc^\prime$, then $\cc^{\prime\prime}$ is a right rejective subcategory of $\cc$}. Since $\epsilon_X$ in the proof of [I5]2.2 was a monomorphism in $\mod\Lambda$, this change has nothing to do with the main results in [I5].)}, where we denote by $\cc_n/[\cc_{n+1}]$ the factor category \ZBA. Dually, a {\it left rejective chain} is defined, and we call a left and right rejective chain a {\it rejective chain}.

For example, the chain $\lat\Lambda_0\subset\cdots\subset\lat\Lambda_{n-1}\subset\lat\Lambda_n$ induced by the overorders $\Lambda_n\subset\Lambda_{n-1}\subset\cdots\subset\Lambda_0$ in \ZAB(2) or (3) is a rejective chain. But the chain $\lat\Lambda_1\subset\cdots\subset\lat\Lambda_{m-1}\subset\lat\Lambda_m$ induced by the overorders $\Lambda_m\subset\Lambda_{m-1}\subset\cdots\subset\Lambda_1$ in \ZAD\ is neither a left nor right rejective chain since the condition $J_{\cc_n/[\cc_{n+1}]}=0$ is not satisfied. But we will see in \ZCDA\ and \ZCGA\ that the category $\lat\Lambda$ for an $R$-order $\Lambda$ and the category $\mod\Lambda$ for an artin algebra $\Lambda$ always have `sufficiently many' right rejective chains. This fact plays a crucial role in this paper.

\vskip1em{\bf\ZBG\ }
The proposition below shows that right rejective chains of additive categories induce heredity chains of quasi-hereditary algebras. For example, the heredity chains of Auslander algebras, which Dlab-Ringel [DR3] gave by using the preprojective partition of Auslander-Smalo [AS1], are induced by right rejective chains.

\vskip1em{\bf Proposition }{\it
Assume that an additive category $\cc$ has an additive generator $M$ and $\Gamma:=\cc(M,M)$ is an artin algebra.

(1) A bijection $\{\cc^\prime:$ right rejective subcategory of $\cc$ with $J_{\cc^\prime}=0\}\rightarrow\{ I:$ heredity ideal of $\Gamma\}$ is given by $\cc^\prime\mapsto I:=[\cc^\prime](M,M)$.

(2) If there exists a right rejective chain $0=\cc_m\subseteq\cc_{m-1}\subseteq\cdots\subseteq\cc_0=\cc$ from $\cc$ to $0$, then $\Gamma$ is a quasi-hereditary algebra with a heredity chain $0=[\cc_m](M,M)\subseteq[\cc_{m-1}](M,M)\subseteq\cdots\subseteq[\cc_0](M,M)=\Gamma$ [I5]. Thus $\gl\Gamma\le 2m-2$ holds [DR1].}

\vskip1em{\sc Proof }
(1) Let $\cc^\prime$ be a right rejective subcategory of $\cc$ with $J_{\cc^\prime}=0$ and $I:=[\cc^\prime](M,M)$. Then $I^2=I$ holds. Take a right $[\cc^\prime]$-approximation $f\in\cc(N,M)$ which is a monomorphism. Then $f$ induces a homomorphism $(\cdot f):\cc(M,N)\to\cc(M,M)$ of $\Gamma$-modules, which is an injection with the image $I$. Thus $I$ is isomorphic to $\cc(M,N)\in\pr\Gamma$. If $g\in\cc(M,M)$ is contained in $IJ_\Gamma I$, then $g$ factors through a morphism in $J_{\cc^\prime}$ which is zero. Thus $IJ_\Gamma I=0$ holds, and $I$ is a heredity ideal of $\Gamma$.

Conversely, let $I$ be a heredity ideal of $\Gamma$. Since $I^2=I$, there exists a subcategory $\cc^\prime$ of $\cc$ such that $I=[\cc^\prime](M,M)$. Since $I\in\pr\Gamma$, there exists an isomorphism $\cc(M,N)\to I$ of $\Gamma$-modules for some $N\in\cc$. The composition $\cc(M,N)\to I\subseteq\cc(M,M)$ is induced by some $f\in\cc(N,M)$. Then $f$ gives a right $[\cc^\prime]$-approximation of $M$ which is a monomorphism. Thus $\cc^\prime$ is a right rejective subcategory of $\cc$.

(2) By \ZBCA(1), $0=\cc_{m-1}/[\cc_{m-1}]\subseteq\cc_{m-2}/[\cc_{m-1}]\subseteq\cdots\subseteq\cc_0/[\cc_{m-1}]=\cc/[\cc_{m-1}]$ is again a right rejective chain. Since $[\cc_{m-1}](M,M)$ is a heredity ideal of $\Gamma$ by (1), the assertion follows inductively.\rule{5pt}{10pt}

\vskip1em{\bf\ZBGA\ Example }Let $\Lambda$ be an artin algebra. Put $\cc_n:=\add\bigoplus_{i=0}^{m-n}\Lambda/J_\Lambda^i$ for $m$ such that $J_\Lambda^m=0$. Then $0=\cc_m\subseteq\cc_{m-1}\subseteq\cdots\subseteq\cc_0$ gives a left rejective chain from $\cc_0$ to $0$. In particular, $\Gamma:=\endm_\Lambda(\bigoplus_{i=0}^{m}\Lambda/J_\Lambda^i)$ is a quasi-hereditary algebra with $\gl\Gamma\le m$, the theorem of Auslander-Dlab-Ringel [A1][DR2].

\vskip1em{\sc Proof }Since any $f\in\cc_0(\Lambda,\Lambda/J_\Lambda^i)$ satisfies $f(J_\Lambda^i)=0$, the natural surjection $p:\Lambda\to\Lambda/J_\Lambda^{m-1}$ gives a left $[\cc_1]$-approximation of $\Lambda$. Thus $\cc_1$ is a left rejective subcategory of $\cc_0$. Moreover, since any $f\in J_{\cc_0}(\Lambda,\Lambda)$ satisfies $f(J_\Lambda^{m-1})=0$, $f$ factors through $p$. Thus $J_{\cc_0/[\cc_1]}=0$ holds. We obtain the assertion by using the dual of \ZBCA(2) repeatedly.\rule{5pt}{10pt}

\vskip1em{\bf\ZBGB\ }Finally, we recall the theorem of Dlab-Ringel which asserts that any artin algebra $\Lambda$ with $\gl\Lambda\le2$ is quasi-hereditary [DR1], which is generalized to orders by K\"onig-Wiedemann [KW]. One can show the proposition below [I9], which is a slightly stronger version of theirs.

\vskip1em{\bf Theorem }{\it Let $\Lambda$ be an artin algebra (resp. $R$-order) with $\gl\Lambda\le2$ and $\cc:=\pr\Lambda$. Then there exists a right rejective chain from $\cc$ to $0$ (resp. from $\cc$ to $\lat\Gamma$ for some maximal overorder $\Gamma$ of $\Lambda$).}

\vskip1em{\bf\ZC\ Representation dimension }

In this section, assume that $\Lambda$ and $\Gamma$ are artin algebras unless stated otherwise. We denote by $\pd_\Gamma X$ the projective dimension of $X\in\mod\Gamma$, and by $0\rightarrow\Gamma\rightarrow I_0(\Gamma)\rightarrow I_1(\Gamma)\rightarrow I_2(\Gamma)\rightarrow\cdots$ the minimal injective resolution of the left $\Gamma$-module $\Gamma$. Then the {\it dominant dimension} of $\Gamma$ is defined by $\domdim\Gamma:=\inf\{i\ge0\ |\ \pd_\Gamma I_i(\Gamma)\neq0\}$ [T]. Let us start with recalling the following classical theorem of Auslander [A1][ARS].

\vskip1em{\bf\ZCA\ (Auslander correspondence) }{\it
There exists a bijection between Morita-equivalence classes of artin algebras $\Lambda$ with $\#\ind(\mod\Lambda)<\infty$ and those of artin algebras $\Gamma$ with $\gl\Gamma\le2$ and $\domdim\Gamma\ge2$. The correspondence is given by $\Lambda\mapsto\Gamma:=\endm_\Lambda(\bigoplus_{X\in\ind(\mod\Lambda)}X)$ and $\Gamma\mapsto\Lambda:=\endm_\Gamma(I_0(\Gamma))$.}

\vskip1em
We call such $\Lambda$ an algebra {\it of finite representation type} and such $\Gamma$ an {\it Auslander algebra}. This quite surprising theorem gives the relationship between a representation theoretic property `finite representation type' and a homological property `Auslander algebra'. It is one of the most important theorem in the representation theory of artin algebras, which leads to later Auslander-Reiten theory [ARS] and is closely related to some duality on Auslander-Gorenstein rings [I6,8].

It would be natural to ask which `representation theoretic' class of algebras corresponds to more general `homological' class than Auslander algebras. Is there some kind of `higher-dimensional' Auslander-Reiten theory which contains usual theory as a `two-dimensional' version? One direction was already given by Auslander about 30 years ago. He introduced the class of algebras with the representation dimension at most $n$, which corresponds to the class of algebras $\Gamma$ with $\gl\Gamma\le n$ and $\domdim\Gamma\ge2$ [A1].

\vskip1em{\bf\ZCB\ }The {\it representation dimension} of an artin algebra $\Lambda$ is defined by
\[\rdim\Lambda:=\inf_{\Gamma}\gl\Gamma,\]
where we consider all artin algebras $\Gamma$ such that $\domdim\Gamma\ge2$ and $\endm_\Gamma(I_0(\Gamma))$ is Morita-equivalent to $\Lambda$. Moreover, Auslander proved basic results below. We often use (1), which is convenient to study $\rdim\Lambda$.

\vskip1em{\bf Theorem }{\it
(1) $\rdim\Lambda=\inf\{\gl\endm_\Lambda(\Lambda\oplus\Lambda^*\oplus N)\ |\ N\in\mod\Lambda\}$.

(2) $\Lambda$ is of finite representation type if and only if $\rdim\Lambda\le2$.}

In the sense of (2), $\rdim\Lambda$ measures how far an artin algebra is from being of finite representation type. Unfortunately, much is unknown about the representation dimension. A reason may be caused by a relationship to some open homological problems, e.g. finitistic dimension conjecture \ZCF. In other words, there might be a possibility that the representation dimension gives some kind of `representation theoretic' approach to such problems, e.g. \ZCFB. By the way, a basic problem \ZCC(1) (conjectured in [X2]) was open until recently, and we will prove it along [I5]. More generally, we will prove the conjecture \ZCC(2) of Ringel-Yamagata, which relates the representation dimension to quasi-hereditary algebras. Putting $M=\Lambda\oplus\Lambda^*$, we can obtain \ZCC(1) from \ZCC(2) and \ZBE\ immediately.

\vskip1em{\bf\ZCC\ Problem }
(1) Is $\rdim\Lambda$ finite for any artin algebra $\Lambda$?

(2) For any artin algebra $\Lambda$ and $M\in\mod\Lambda$, is there $N\in\mod\Lambda$ such that $\endm_\Lambda(M\oplus N)$ is a quasi-hereditary algebra?

\vskip1em{\bf\ZCD\ }
To prove \ZCC(2), we only have to construct a subcategory $\cc$ of $\mod\Lambda$ such that $M\in\cc$ and there exists a right rejective chain from $\cc$ to $0$ by \ZBG. This is done by \ZCDA\ below [I5]. Thus we obtain an answer \ZCDB\ to \ZCC.

\vskip1em{\bf\ZCDA\ Theorem }{\it Let $\Lambda$ be an artin algebra and $M\in\mod\Lambda$. Put $M_0:=M$, $M_{n+1}:=M_nJ_{\endm_\Lambda(M_n)}\subsetneq M_n$ and take large $m$ such that $M_m=0$. Then $\cc_n:=\add\bigoplus_{l=n}^{m-1}M_l$ gives a right rejective chain $0=\cc_m\subseteq\cc_{m-1}\subseteq\cdots\subseteq\cc_0$ from $\cc_0$ to $0$. Thus $\endm_\Lambda(\bigoplus_{l=0}^{m-1}M_l)$ is a quasi-hereditary algebra.}

\vskip1em{\sc Proof }
(i) Note that there exists a surjection $f_{n,l}\in\hom_\Lambda(\bigoplus M_n,M_l)$ for any $n<l$.

(ii) Define a functor $\fff_n:\mod\Lambda\rightarrow\mod\Lambda$ by
\[\fff_n(X):=\sum_{Y\in\cc_n,\ f\in J_{\mod\Lambda}(Y,X)}f(Y).\]
Then a natural transformation $\epsilon:\fff_n\rightarrow 1$ is defined by the inclusion $\epsilon_X:\fff_n(X)\rightarrow X$. By (i), $\fff_n(M_n)=M_nJ_{\endm_\Lambda(M_n)}=M_{n+1}\in\cc_{n+1}$ holds. Thus $J_{\cc_n}(\ ,X)=[\cc_{n+1}](\ ,X)=\cc_n(\ ,\fff_n(X))\epsilon_X$ holds on $\cc_n$ for any $X\in\ind\cc_n-\ind\cc_{n+1}$.

(iii) Fix $X\in\ind\cc_n$. Put $Y:=\fff_n(X)$ and $g:=\epsilon_X$ if $X\notin\cc_{n+1}$, and $Y:=X$ and $g:=1_X$ if $X\in\cc_{n+1}$. By (ii), $g\in\cc_n(Y,X)$ gives a right $[\cc_{n+1}]$-approximation of $X$ which is a monomorphism in $\mod\Lambda$. Thus $\cc_{n+1}$ is a right rejective subcategory of $\cc_n$, and $J_{\cc_n/[\cc_{n+1}]}=0$ holds by (ii). We obtain the assertion by using \ZBCA(2) repeatedly.\rule{5pt}{10pt}

\vskip1em{\bf\ZCDB\ Corollary }{\it Let $\Lambda$ be an artin algebra and $M\in\mod\Lambda$. Then there exists $N\in\mod\Lambda$ such that $\endm_\Lambda(M\oplus N)$ is a quasi-hereditary algebra. Thus $\rdim\Lambda<\infty$.}

\vskip.5em{\bf\ZCDC\ }We will give one more application of \ZCDA. It is not difficult to show that, if a subcategory $\cc$ of $\mod\Lambda$ with $\#\ind\cc<\infty$ is closed under submodules, then $\gl\endm_\Lambda(\bigoplus_{X\in\ind\cc}X)\le2$ holds, and there exists a right rejective chain from $\cc$ to $0$ by \ZBGB. By \ZCDA, we can show its funny generalization below [I9]. For example, $\cc^{(n)}:=\add\bigoplus_{X\in\ind(\mod\Lambda),\ \length_\Lambda X\le n}X$ satisfies the assumption for any $n\in\nnn$.

\vskip1em{\bf Theorem }{\it Let $\Lambda$ be an artin algebra and $\cc$ a subcategory of $\mod\Lambda$ such that $\#\ind\cc<\infty$. Assume that any submodule of any $X\in\ind\cc$ is contained in $\cc$. Then there exists a right rejective chain from $\cc$ to $0$. Thus $\endm_\Lambda(\bigoplus_{X\in\ind\cc}X)$ is a quasi-hereditary algebra.}

\vskip1em{\bf\ZCE\ }Here we collect some results on the representation dimension. See also [H][X3].

(1) If $\gl\Lambda\le1$ or $J_\Lambda^2=0$ holds, then $\rdim\Lambda\le3$ [A1].

(2) If $\Lambda$ is selfinjective, then $\rdim\Lambda\le($the Loewy length of $\Lambda)$ [A1] by \ZBGA.

(3) If $\Lambda$ and $\Gamma$ are algebras over a perfect field $k$, then $\rdim(\Lambda\otimes_k\Gamma)\le\rdim\Lambda+\rdim\Gamma$ [X1]. In particular, $\rdim\left(\def\arraystretch{.5}\begin{array}{cc}{\scriptstyle\Lambda}&{\scriptstyle\Lambda}\\ {\scriptstyle0}&{\scriptstyle\Lambda}\end{array}\right)\le\rdim\Lambda+2$ [FGR].

(4) Stable equivalences of Morita type preserve the representation dimension [X2].

(5) If $\Lambda$ is a special biserial algebra, then $\rdim\Lambda\le3$ [EHIS].

\vskip1em{\bf\ZCF\ Problem }
What is the set $\{\rdim\Lambda\ |\ \Lambda$ is an artin algebra$\}$? This is unknown. It has been shown that some classes of artin algebras $\Lambda$ satisfy $\rdim\Lambda\le 3$ (\ZCE). It is curious that a question of Auslander [A1], whether any artin algebra $\Lambda$ satisfies $\rdim\Lambda\le3$ or not, is still unknown.\footnote{At a conference in November 2002, R. Rouquier announced a proof that the exterior algebra of a 3-dimensional vector space has representation dimension 4.} Nevertheless we will see that the condition $\rdim\Lambda\le3$ is very interesting.

Now we define the {\it finitistic dimension} [B] of $\Lambda$ by
\[\fdim\Lambda:=\sup\{\pd_\Lambda X\ |\ X\in\mod\Lambda,\ \pd_\Lambda X<\infty\}.\]
The {\it finitistic dimension conjecture}, which asserts that any artin algebra $\Lambda$ satisfies $\fdim\Lambda<\infty$, is an old open problem [Z]. This is a rather strong conjecture among many open homological problems, e.g. Nakayama conjecture. Recently, Igusa-Todorov proved the following interesting result [IT].

\vskip1em{\bf\ZCFA\ Theorem }{\it
Let $\Gamma$ be an artin algebra with $\gl\Gamma\le3$ and $P\in\pr\Gamma$. Then $\Lambda:=\endm_\Gamma(P)$ satisfies $\fdim\Lambda<\infty$.}

\vskip1em
As an easy conclusion, $\rdim\Lambda\le3$ implies $\fdim\Lambda<\infty$. (Take $N\in\mod\Lambda$ such that $\Gamma:=\endm_\Lambda(\Lambda\oplus\Lambda^*\oplus N)$ satisfies $\gl\Gamma\le3$. Then $P:=\hom_\Lambda(\Lambda\oplus\Lambda^*\oplus N,\Lambda)\in\pr\Gamma$ satisfies $\Lambda=\endm_\Gamma(P)$.) Thus we obtain the corollary [EHIS] below by \ZCE(5). See [X4] for more results along this approach.

\vskip1em{\bf\ZCFB\ Corollary }{\it $\fdim\Lambda<\infty$ holds for any special biserial algebra $\Lambda$.}

\vskip1em{\bf\ZCG\ Representation dimension of orders }

Let $R$ be a complete discrete valuation ring and $\Lambda$ an $R$-order in a semisimple algebra $A=\Lambda\otimes_RK$. The one-dimensional analogy of Auslander correspondence was given by Auslander-Roggenkamp [AR]. Now we define the {\it representation dimension} of $\Lambda$ by
\[\rdim\Lambda:=\inf\{\gl\endm_\Lambda(\Lambda\oplus\Lambda^*\oplus N)\ |\ N\in\lat\Lambda\}.\]
Then $\Lambda$ is of finite representation type (i.e. $\#\ind(\lat\Lambda)<\infty$) if and only if $\rdim\Lambda\le2$. Moreover, we can show the following analogy of \ZCD\ [I9].

\vskip1em{\bf\ZCGA\ Theorem }{\it Let $\Lambda$ be an $R$-order in a semisimple algebra $A$ and $M\in\lat\Lambda$.

(1) Put $M_0:=M$ and $M_{n+1}:=M_nJ_{\endm_\Lambda(M_n)}\subsetneq M_n$. Then there exists $m\ge0$ such that $\cc_n:=\add\bigoplus_{l=n}^{m-1}M_l$ gives a right rejective chain $\cc_m\subseteq\cc_{m-1}\subseteq\cdots\subseteq\cc_0$ from $\cc_0$ to $\cc_m=\lat\Gamma$ for a hereditary overring $\Gamma$ of $\Lambda$. Thus $\endm_\Lambda(\bigoplus_{l=0}^{m-1}M_l)$ is a quasi-hereditary order.

(2) There exists $N\in\mod\Lambda$ such that $\endm_\Lambda(M\oplus N)$ is a quasi-hereditary order. Thus $\rdim\Lambda<\infty$.}

\vskip1em{\bf\ZCGB\ Question }
It would be natural to define a higher-dimensional analogy of the representation dimension as follows: Let $R$ be a complete regular local ring, $\Lambda$ an $R$-order and $\lat\Lambda$ the category of $\Lambda$-lattices (\S1). Define the {\it representation dimension} of $\Lambda$ by $\rdim\Lambda=\inf\{\gl\endm_\Lambda(\Lambda\oplus\Lambda^*\oplus N)\ |\ N\in\lat\Lambda\}$. Is $\rdim\Lambda$ finite?
What is the possible value of $\rdim\Lambda$?

\vskip1em{\bf\ZD\ Solomon zeta function }

Let $R$ be the ring $\zzz$ of integers or its $p$-adic completion $\zzz_p$, and $K$ its quotient field. For an $R$-order $\Lambda$ in a semisimple $K$-algebra $A$, its {\it Solomon zeta function} is defined by $\zeta_\Lambda(s):=\sum_{L}(\Lambda:L)^{-s}$ where $L$ is a left ideal of $\Lambda$ such that $(\Lambda:L)<\infty$ and $s$ is a complex variable [S1]. Then $\zeta_\Lambda$ converges in the half-plane $\{ s\in\ccc\ |\ \real(s)>\dim_KA\}$, and it can be shown that $\zeta_\Lambda$ admits analytic continuation to a meromorphic function of $s$ [S1]. When $\Lambda$ is the ring of integers in an algebraic number field $A$, then $\zeta_\Lambda$ is the usual Dedekind zeta function of $A$. For the case $R=\zzz$, the {\it Euler product formula} $\zeta_\Lambda=\prod_{p:{\rm prime}}\zeta_{\Lambda_p}$ holds for $\Lambda_p:=\zzz_p\otimes_{\zzz}\Lambda$ [S1]. If $\Lambda$ is a Gorenstein order (i.e. $\Lambda\in\rin\Lambda$) with a maximal overorder $\Gamma$, then the {\it functional equation} $\zeta_\Lambda(s)/\zeta_\Lambda(1-s)=(\Gamma:\Lambda)^{1-2s}\zeta_\Gamma(s)/\zeta_\Gamma(1-s)$ holds [G]. Later, Bushnell-Reiner applied Solomon zeta functions to obtain the prime ideal theorem and the asymptotic distribution formula of ideals for general orders [BR1,2,3,4], which were well-known for the rings of integers in algebraic number fields [L].

For example, $\zeta_{\ma_n(\zzz)}(s)$ equals $\prod_{i=0}^{n-1}\zeta_{\zzz}(ns-i)$. More generally, for a maximal overorder $\Gamma$ of $\Lambda$, Hey gave an explicit description of $\zeta_\Gamma$, which is a product of Dedekind zeta functions and correction factors [BR2]\S2. Since $\{ p:{\rm prime}\ |\ \Lambda_p\neq\Gamma_p\}$ is a finite set, the difference between $\zeta_\Lambda$ and $\zeta_\Gamma$ appears at only finitely many primes. To study it, we will consider the case $R=\zzz_p$ in the rest of this section.

Thus, in the rest of this section, let $R$ be a complete discrete valuation ring with the residue field $k$ and the quotient field $K$. Assume that $k$ is a finite field with $p$ elements. It will be natural to define the zeta function not only for ideals but also for modules.

\vskip1em{\bf\ZDA\ Definition }
In the rest, let $\Lambda$ be an $R$-order in a semisimple algebra $A=\Lambda\otimes_RK$. We denote by $\widetilde{(\ )}$ the functor $(\ )\otimes_RK:\lat\Lambda\rightarrow\mod A$. For $V\in\mod A$, we denote by $\ll_\Lambda(V)$ the partially ordered set of {\it full} $\Lambda$-lattices in $V$, which are $\Lambda$-lattices $L$ such that $L\subset V$ and $KL=V$. Then the set of isoclasses $\overline{\ll}_\Lambda(V):=\ll_\Lambda(V)/\simeq$ is a finite set by Jordan-Zassenhaus Theorem [CR]. For $L,M\in\overline{\ll}_\Lambda(V)$, Solomon [S1] studied
\begin{eqnarray*}
\mbox{a {\it partial zeta function}}&&{\rm Z}(L,M;s):=\sum_{N\subseteq L,\ N\simeq M}(L:N)^{-s}\\
\mbox{and the {\it zeta matrix} of size $\#\overline{\ll}_\Lambda(V)$}&&{\bf Z}_{\Lambda}(V;s):=({\rm Z}(L,M;s))_{L,M\in\overline{\ll}_\Lambda(V)}.
\end{eqnarray*}
He proved that ${\bf Z}_{\Lambda}(V;s)$ has an inverse matrix in $\ma_n(\zzz[p^{-s}])$ for $n:=\#\overline{\ll}_\Lambda(V)$ by a combinatorial argument (M\"obius inversion) [S1]. In particular, ${\rm Z}(L,M;s)$ and $\zeta_\Lambda(s)$ are rational functions of $p^{-s}$. Moreover, he gave the following conjectures in [S2].

\vskip1em{\bf\ZDB\ Problem }
(1) ${\rm Z}(L,M;s)/\det{\bf Z}_{\Gamma}(V;s)\in\zzz[p^{-s}]$ for a maximal order $\Gamma$ in $A$.

(2) $\det{\bf Z}_{\Lambda}(V;s)$ is the finite product $\prod_{i}(1-p^{a_i-b_is})^{-1}$ with some $a_i\in\nnn_{\ge0}$ and $b_i\in\nnn_{>0}$.

\vskip1em
Solomon's first conjecture (1) was proved in [BR1] by using zeta integrals. However Solomon's second conjecture (2) was open until recently, although a special case when $\Lambda$ is hereditary was proved by Denert [D]. In the rest of this section, we will give an explicit description of $\det{\bf Z}_{\Lambda}(V;s)$ for general $\Lambda$ [I4], which implies the second conjecture.

\vskip1em{\bf\ZDC\ Example }
(1) Put $\Lambda:=\left(\def\arraystretch{.5}\begin{array}{ccccc}{\scriptstyle R}&{\scriptstyle R}&{\scriptstyle R}&{\scriptstyle \cdots}&{\scriptstyle R}\\ {\scriptstyle J_R}&{\scriptstyle R}&{\scriptstyle R}&{\scriptstyle \cdots}&{\scriptstyle R}\\ {\scriptstyle J_R}&{\scriptstyle J_R}&{\scriptstyle R}&{\scriptstyle \cdots}&{\scriptstyle R}\\ {\scriptstyle \cdots}&{\scriptstyle \cdots}&{\scriptstyle \cdots}&{\scriptstyle \cdots}&{\scriptstyle \cdots}\\ {\scriptstyle J_R}&{\scriptstyle J_R}&{\scriptstyle J_R}&{\scriptstyle \cdots}&{\scriptstyle R} \end{array}\right)\subset A:=\ma_n(K)$ and $V:=\ma_{n,1}(K)$. Then $\overline{\ll}_\Lambda(V)=\{ P_i\}_{1\le i\le n}$ holds for the $i$-th row $P_i$ of $\Lambda$. Putting $T:=p^{-s}$, we obtain \[{\bf Z}_\Lambda(V;s)=(1-T^n)^{-1}\left(\def\arraystretch{.5}\begin{array}{ccccc}{\scriptstyle 1}&{\scriptstyle T^{n-1}}&{\scriptstyle T^{n-2}}&{\scriptstyle \cdots}&{\scriptstyle T}\\ {\scriptstyle T}&{\scriptstyle 1}&{\scriptstyle T^{n-1}}&{\scriptstyle \cdots}&{\scriptstyle T^2}\\ {\scriptstyle T^2}&{\scriptstyle T}&{\scriptstyle 1}&{\scriptstyle \cdots}&{\scriptstyle T^3}\\ {\scriptstyle \cdots}&{\scriptstyle \cdots}&{\scriptstyle \cdots}&{\scriptstyle \cdots}&{\scriptstyle \cdots}\\ {\scriptstyle T^{n-1}}&{\scriptstyle T^{n-2}}&{\scriptstyle T^{n-3}}&{\scriptstyle \cdots}&{\scriptstyle 1}\end{array}\right).\]

Thus $\det{\bf Z}_\Lambda(V;s)=(1-T^n)^{-1}$ holds.

(2) Let $\Lambda_n$ be the order in \ZAB(2) and $T:=p^{-s}$. Then $\overline{\ll}_{\Lambda_n}(A)=\{\Lambda_i\}_{0\le i\le n}$ holds. For example, ${\bf Z}_{\Lambda_3}(A;s)$ is
{\tiny \[(1-T)^{-2}\left(\begin{array}{rrrr}
{\scriptscriptstyle \begin{array}{lr}1-2T+(p+1)T^2-2pT^3\\
\ \ \ \ \ +(p^2+p)T^4-2p^2T^5+p^3T^6\end{array}}&{\scriptscriptstyle T-2T^2+(p+1)T^3-2pT^4+p^2T^5}&{\scriptscriptstyle T^2-2T^3+pT^4}&{\scriptscriptstyle T^3}\\
{\scriptscriptstyle pT-2pT^2+(p^2+p)T^3-2p^2T^4+p^3T^5}&{\scriptscriptstyle 1-2T+(p+1)T^2-2pT^3+p^2T^4}&{\scriptscriptstyle T-2T^2+pT^3}&{\scriptscriptstyle T^2}\\
{\scriptscriptstyle p^2T^2-2p^2T^3+p^3T^4}&{\scriptscriptstyle pT-2pT^2+p^2T^3}&{\scriptscriptstyle 1-2T+pT^2}&{\scriptscriptstyle T}\\
{\scriptscriptstyle (p^3-p^2)T^3}&{\scriptscriptstyle (p^2-p)T^2}&{\scriptscriptstyle (p-1)T}&{\scriptscriptstyle 1}\end{array}\right).\]}

By an elementary transformation of rows, ${\bf Z}_{\Lambda_3}(A;s)$ is chaged into
{\tiny \[(1-T)^{-2}\left(\begin{array}{rrrr}
{\scriptstyle 1-2T+T^2}&{\scriptstyle 0}&{\scriptstyle 0}&{\scriptstyle 0}\\
{\scriptstyle pT-2pT^2+(p^2+p)T^3-2p^2T^4+p^3T^5}&{\scriptstyle 1-2T+(p+1)T^2-2pT^3+p^2T^4}&{\scriptstyle T-2T^2+pT^3}&{\scriptstyle T^2}\\
{\scriptstyle p^2T^2-2p^2T^3+p^3T^4}&{\scriptstyle pT-2pT^2+p^2T^3}&{\scriptstyle 1-2T+pT^2}&{\scriptstyle T}\\
{\scriptstyle (p^3-p^2)T^3}&{\scriptstyle (p^2-p)T^2}&{\scriptstyle (p-1)T}&{\scriptstyle 1}\end{array}\right),\]}
which is equal to $\left(\def\arraystretch{.5}\begin{array}{cc}
{\scriptstyle 1}&{\scriptstyle O}\\
{*}&{\bf Z}_{\Lambda_2}(A;s)\end{array}\right)$. Thus $\det{\bf Z}_{\Lambda_3}(A;s)=\det{\bf Z}_{\Lambda_2}(A;s)$ holds. By a quite similar argument, we obtain $\det{\bf Z}_{\Lambda_n}(A;s)=\det{\bf Z}_{\Lambda_0}(V;s)=(1-T)^{-2}$.

(3) Let $\Lambda_n$ be the order in \ZAB(3) and $T:=p^{-s}$. Then $\overline{\ll}_{\Lambda_n}(A)=\{\Lambda_i\}_{0\le i\le n}$ holds. For example, ${\bf Z}_{\Lambda_3}(A;s)$ is
{\tiny \[(1-T)^{-1}\left(\begin{array}{rrrr}
{\scriptstyle 1-T+pT^2-pT^3+p^2T^4-p^2T^5+p^3T^6}&{\scriptstyle T-T^2+pT^3-pT^4+p^2T^5}&{\scriptstyle T^2-T^3+pT^4}&{\scriptstyle T^3}\\
{\scriptstyle pT-pT^2+p^2T^3-p^2T^4+p^3T^5}&{\scriptstyle 1-T+pT^2-pT^3+p^2T^4}&{\scriptstyle T-T^2+pT^3}&{\scriptstyle T^2}\\
{\scriptstyle p^2T^2-p^2T^3+p^3T^4}&{\scriptstyle pT-pT^2+p^2T^3}&{\scriptstyle 1-T+pT^2}&{\scriptstyle T}\\
{\scriptstyle p^3T^3}&{\scriptstyle p^2T^2}&{\scriptstyle pT}&{\scriptstyle 1}\end{array}\right).\]}

By an elementary transformation of rows, ${\bf Z}_{\Lambda_3}(A;s)$ is chaged into
{\tiny \[(1-T)^{-1}\left(\begin{array}{rrrr}
{\scriptstyle 1-T}&{\scriptstyle 0}&{\scriptstyle 0}&{\scriptstyle 0}\\
{\scriptstyle pT-pT^2+p^2T^3-p^2T^4+p^3T^5}&{\scriptstyle 1-T+pT^2-pT^3+p^2T^4}&{\scriptstyle T-T^2+pT^3}&{\scriptstyle T^2}\\
{\scriptstyle p^2T^2-p^2T^3+p^3T^4}&{\scriptstyle pT-pT^2+p^2T^3}&{\scriptstyle 1-T+pT^2}&{\scriptstyle T}\\
{\scriptstyle p^3T^3}&{\scriptstyle p^2T^2}&{\scriptstyle pT}&{\scriptstyle 1}\end{array}\right),\]}
which is equal to $\left(\def\arraystretch{.5}\begin{array}{cc}
{\scriptstyle 1}&{\scriptstyle O}\\
{*}&{\bf Z}_{\Lambda_2}(A;s)\end{array}\right)$. Thus $\det{\bf Z}_{\Lambda_3}(A;s)=\det{\bf Z}_{\Lambda_2}(A;s)$ holds. By a quite similar argument, we obtain $\det{\bf Z}_{\Lambda_n}(A;s)=\det{\bf Z}_{\Lambda_0}(V;s)=(1-T)^{-1}$.

\vskip1em
As the above examples (2) and (3) show, we can often use a `good' sequence of overorders of $\Lambda$ to calculate $\det{\bf Z}_\Lambda(A;s)$. In general, there does not necessarily exist such a `good' sequence. But, we will see below that we can always use a right rejective chain, which can be regarded as a generalization of a `good' sequence of overorders (\ZBF).

\vskip1em{\bf\ZDD\ }Let us follow the approach in [I4]. Our first step is to define the zeta matrix for any subcategory $\cc$ (\S\ZB) of $\lat\Lambda$. We denote by $\ll_{\cc}(V)$ (resp. $\overline{\ll}_{\cc}(V)$) the subset of $\ll_\Lambda(V)$ (resp. $\overline{\ll}_\Lambda(V)$) consisting of objects in $\cc$. We will study
\begin{eqnarray*}\mbox{the {\it zeta matrix}}&&{\bf Z}_{\cc}(V;s):=({\rm Z}(L,M;s))_{L,M\in\overline{\ll}_{\cc}(V)}
\end{eqnarray*}
which has the size $\#\overline{\ll}_{\cc}(V)$. For a general subcategory $\cc$, $\det{\bf Z}_{\cc}(V;s)$ is far from the form in \ZDB(2). But we will see in \ZDF(2) below that if there exists a right rejective chain (\ZBF) from $\cc$ to $\lat\Gamma$ for some maximal overring $\Gamma$ of $\Lambda$, then we can give an explicit description of $\det{\bf Z}_{\cc}(V;s)$, which has the form in \ZDB(2). On the other hand, the subcategory $\cc_V$ for any $V\in\mod A$ in \ZDE\ below has a such right rejective chain. Consequently, ${\bf Z}_{\Lambda}(V;s)$ itself also has the form in \ZDB(2) since ${\bf Z}_{\Lambda}(V;s)={\bf Z}_{\cc_V}(V;s)$ holds.

\vskip1em{\bf\ZDE\ Theorem }{\it (cf. \ZCDC) Let $\Lambda$ be an $R$-order in a semisimple algebra $A$ and $V\in\mod A$. Put $\cc_V:=\add\bigoplus_{X\in\ind(\lat\Lambda),\ X\subset V}X$. Then there exists a right rejective chain from $\cc_V$ to $\lat\Gamma$ for a maximal overring $\Gamma$ of $\Lambda$, and $\endm_\Lambda(\bigoplus_{X\in\ind\cc_V}X)$ is a quasi-hereditary order.}

\vskip1em{\sc Proof }
We give a different proof from [I4]2.6 as an application of \ZCGA(1). Apply the construction of \ZCGA(1) to $M:=\bigoplus_{X\in\ind\cc_V}X$. Since $M_{n+1}$ is given by $M_{n+1}=M_nJ_{\endm_\Lambda(M_n)}$, we obtain $M_n\in\cc_V$ for any $n$ inductively. Thus $\cc_V=\add\bigoplus_{l=0}^{m-1}M_l$ holds. The assertion follows from \ZCGA(1) immediately.
\rule{5pt}{10pt}

\vskip1em{\bf\ZDF\ }Let us state the main theorem. Put $A=\prod_{j=1}^rA_j$ for simple algebras $A_j$. Let $e_j$ be the identity of $A_j$, $\Gamma_j$ a maximal overorder of $e_j\Lambda$ in $A_j$, $\Gamma:=\prod_{j=1}^r\Gamma_j$, $S_j$ a simple $A_j$-module, and $G_j$ a simple $\Gamma_j$-module. Then $(S_j)_{1\le j\le r}$ (resp. $(G_j)_{1\le j\le r}$) gives the set of isoclasses of simple $A$-modules (resp. simple $\Gamma$-modules). For $X\in\mod A$, we denote by $l_j(X)$ the multiplicity of $S_j$ as a composition factor of $X$. Put $q_j:=\#\endm_\Gamma(G_j)=p^{\dim_k\endm_\Gamma(G_j)}$. Then $q_j$ is independent of a choice of $\Gamma_j$ since any maximal order in $A_j$ is conjugate to $\Gamma_j$.

\vskip1em{\bf Theorem }{\it
Let $\Lambda$ be an $R$-order in a semisimple algebra $A$, $V\in\mod A$ and
$V_j:=V/S_j^{l_j(V)}$ ($1\le j\le r$).

(1) The equation below holds. Thus Solomon's second conjecture holds.
\begin{eqnarray*}\det{\bf Z}_{\Lambda}(V;s)=\prod_{j=1}^r\prod_{n=0}^{l_j(V)-1}(1-q_j^{n-l_j(A)s})^{-\#\overline{\ll}_\Lambda(S_j^n\oplus V_j)}.\end{eqnarray*}

(2) More generally, let $\cc$ be a subcategory of $\lat\Lambda$ with a right rejective chain from $\cc$ to $\lat\Gamma$ for a maximal overring $\Gamma$ of $\Lambda$. Then the equation below holds.
\begin{eqnarray*}\det{\bf Z}_{\cc}(V;s)=\prod_{j=1}^r\prod_{n=0}^{l_j(V)-1}(1-q_j^{n-l_j(A)s})^{-\#\overline{\ll}_{\cc}(S_j^n\oplus V_j)}.\end{eqnarray*}}

\vskip0em{\bf\ZDG\ }
In the rest, we will give a proof of \ZDF(2) {\it under the assumption $r=1$, namely $A$ is simple}. Thanks to this assumption, our calculation becomes much simpler than the general one in [I4]. Put $S:=S_1$, $G:=G_1$, $l:=l_1$ and $q:=q_1$. Moreover, put $(L:M):=(L:L\cap M)\cdot(M:L\cap M)^{-1}$ for $V\in\mod K$ and $L,M\in\ll_R(V)$ for simplicity. This symbol is skew-symmetric, and satisfies $(L:M)\cdot(M:N)=(L:N)$.

\vskip1em{\bf\ZDGA\ }For any $N\in\lat\Lambda$, there exists a map $b^N:\overline{\ll}_\Lambda(V)\rightarrow\rrr_{>0}$ such that
\begin{eqnarray*}b^N_L\cdot(b^N_M)^{-1}=(\hom_\Lambda(N,L):\hom_\Lambda(N,M))\cdot(L:M)^{-l(\widetilde{N})/l(A)}\end{eqnarray*}
holds for any $L,M\in\ll_\Lambda(V)$.

\vskip1em{\sc Proof }
Fix $X\in\ll_\Lambda(V)$. For $L\in\ll_\Lambda(V)$, put
\[b^N_L:=(\hom_\Lambda(N,L):\hom_\Lambda(N,X))\cdot(L:X)^{-l(\widetilde{N})/l(A)}.\]
Then it is not difficult to show that $L\simeq M$ implies $b^N_L=b^N_M$.\rule{5pt}{10pt}

\vskip1em{\bf\ZDGB\ }The following key lemma is a generalization of calculations in \ZDC(2) and (3).

\vskip1em{\bf Lemma }{\it
Assume that $\cc^\prime\subset\cc$ is a consecutive two terms in a right rejective chain such that $\ind\cc-\ind\cc^\prime=\{ X\}$. For any $V\in\mod A$ such that $X\subset V$, there exists a diagonal matrix ${\bf B}$ and a nilpotent matrix ${\bf C}$ such that the following equation holds.
\[(1-{\bf C})\cdot{\bf Z}_{\cc}(V;s)=\left(\begin{array}{cc}
{\bf B}\cdot{\bf Z}_{\cc}(V/\widetilde{X};s-l(\widetilde{X})/l(A))\cdot{\bf B}^{-1}&O\\
{*}&{\bf Z}_{\cc^\prime}(V;s)
\end{array}\right)\]

${\bf B}$ and ${\bf C}$ are given by ${\bf B}_{L,L}=b^{X}_L$ and ${\bf C}_{X\oplus L,Y\oplus L}=(X:Y)^{-s}$ for any $L\in\overline{\ll}_{\cc}(V/\widetilde{X})$ and other entries are $0$.}

\vskip1em{\bf\ZDGC\ Sketch of the proof of \ZDF(2)}
Take a right rejective chain $\lat\Gamma=\cc_{m}\subseteq\cdots\subseteq\cc_1\subseteq\cc_0=\cc$. By an easy argument chopping the chain, we can assume $\#(\ind\cc_n-\ind\cc_{n+1})=1$ for any $n$. For the case $m=0$ i.e. $\cc=\lat\Gamma$, we can show \ZDF(2) easily (e.g. the simplest case of Hey's formula [BR1]). Now we assume that \ZDF(2) holds for $\cc^\prime:=\cc_1$. We will show that \ZDF(2) holds for $\cc$ as well by using the induction on the length of $V$. We will apply \ZDGB, where $\det(1-{\bf C})=1$ holds.
\begin{eqnarray*}
&&\det{\bf Z}_{\cc}(V;s)\stackrel{\ZDGB}{=}
\det{\bf Z}_{\cc}(V/\widetilde{X};s-l(\widetilde{X})/l(A))\cdot\det{\bf Z}_{\cc^\prime}(V;s)\\
&=&\prod_{n=0}^{l(V/\widetilde{X})-1}(1-q^{n-l(A)(s-l(\widetilde{X})/l(A))})^{-\#\overline{\ll}_{\cc}(S^n)}\cdot\prod_{n=0}^{l(V)-1}(1-q^{n-l(A)s})^{-\#\overline{\ll}_{\cc^\prime}(S^n)}\\
&=&\prod_{n=0}^{l(V)-1}(1-q^{n-l(A)s})^{-\#\overline{\ll}_{\cc}(S^{n-l(\widetilde{X})})-\#\overline{\ll}_{\cc^\prime}(S^n)}\\
&=&\prod_{n=0}^{l(V)-1}(1-q^{n-l(A)s})^{-\#\overline{\ll}_{\cc}(S^n)},\end{eqnarray*}
where we used Krull-Schmidt theorem to obtain $\#\overline{\ll}_{\cc}(S^{n-l(\widetilde{X})})+\#\overline{\ll}_{\cc^\prime}(S^n)=\#\overline{\ll}_{\cc}(S^n)$.\rule{5pt}{10pt}

\vskip1em{\bf\ZE\ Solomon zeta functions and Ringel-Hall algebras }

In this final section, we relates Solomon zeta functions to Ringel-Hall algebras. Again let $R$ be a complete discrete valuation ring with the residue field $k$ and the quotient field $K$, and we assume that $k$ is a finite field with $p$ elements. Let $\Lambda$ be an $R$-order in a semisimple algebra $A=\Lambda\otimes_RK$. Then any $R$-order $\Lambda$ is finitary in the sense of Ringel [R1]. For $\Lambda$-modules $X,Y,Z$, we denote by $\ff^Y_{XZ}$ the set of submodules $W$ of $Y$ such that $W\simeq Z$ and $Y/W\simeq X$. Put $F^Y_{XZ}:=\#\ff^Y_{XZ}\in\nnn_{\ge0}\cup\{\infty\}$. We denote by $\fin\Lambda$ the set of isoclasses of $\Lambda$-modules of finite length. The {\it Hall algebra} $\hall(\Lambda)$ is defined as a free $\zzz$-module with the basis $( u_{X})_{X\in\fin\Lambda}$ and the multiplication $ u_{X} u_{Z}:=\sum_{Y\in\fin\Lambda}F^{Y}_{XZ} u_{Y}$. Then $\hall(\Lambda)$ forms an associative ring with the identity $ u_{0}$ [R1].

Now we will construct a family of $\hall(\Lambda)$-modules by a Hall algebra approach. For $V\in\mod A$, let $\hallm_V(\Lambda)$ be a free $\zzz$-module with the basis $( u_{L})_{L\in\overline{\ll}_\Lambda(V)}$ (\ZDA). Define the action of $ u_{X}\in\hall(\Lambda)$ on $\hallm_V(\Lambda)$ by $ u_{X} u_{M}:=\sum_{L\in\overline{\ll}_\Lambda(V)}F^{L}_{XM} u_{L}$ for any $M\in\overline{\ll}_\Lambda(V)$. Then $\hallm_V(\Lambda)$ forms an $\hall(\Lambda)$-module, which we call the {\it Hall module}.

\vskip1em{\bf\ZEA\ }
Put $T:=p^{-s}$, $\widehat{\hall}(\Lambda):=\plim_{n\ge0}\hall(\Lambda)\otimes_{\zzz}\zzz[T]/(T^n)$ and $\widehat{\hallm}_V(\Lambda):=\hallm_V(\Lambda)\otimes_{\zzz}\zzz[[T]]$. Then $\widehat{\hall}(\Lambda)$ acts on $\widehat{\hallm}_V(\Lambda)$ naturally. The proposition below shows a relationship between the zeta matrix and Hall algebras, where $z$ is regarded as an analogy of Hecke series [M].

\vskip1em{\bf Proposition }{\it Put $z:=\sum_{X\in\fin\Lambda} u_{X}\otimes (\# X)^{-s}\in\widehat{\hall}(\Lambda)$. Then the action of $z$ on the finite rank $\zzz[[T]]$-module $\widehat{\hallm}_V(\Lambda)$ is given by the zeta matrix ${\bf Z}_{\Lambda}(V;s)$.}

\vskip1em{\sc Proof }For $M\in\overline{\ll}_\Lambda(V)$, the following equation shows the assertion.
\begin{eqnarray*}
z\cdot u_{M}&=&\sum_{X\in\fin\Lambda} u_{X} u_{M}\otimes (\# X)^{-s}\\
&=&\sum_{L\in\overline{\ll}_\Lambda(V)}u_{L}\otimes\sum_{X\in\fin\Lambda}F^L_{XM}(\# X)^{-s}\\
&=&\sum_{L\in\overline{\ll}_\Lambda(V)}u_{L}\otimes{\rm Z}(L,M;s)\ \ \rule{5pt}{10pt}
\end{eqnarray*}

\vskip1em{\bf\ZEB\ Hall algebras of hereditary orders and $U(\dn{gl}_n)$.\ }

In the rest, we follow [I7]. Recall that an order $\Lambda$ is called {\it hereditary} if $\gl\Lambda=1$. Let $\Delta$ be the cyclic quiver with $n$-vertices, whose set $\Delta_0$ of vertices are identified with $\zzz/n\zzz$. We say that an order $\Lambda$ is of {\it type $\Delta$} if $\Lambda$ is Morita-equivalent to $\tri_n(\Omega)=\{ (x_{ij})\in\ma_n(\Omega)\ |\ x_{ij}\in J_\Omega$ for any $i<j\}$ for a maximal order $\Omega$ in a division algebra. Then put $q_\Lambda:=\#(\Omega/J_\Omega)$. It is well-known that any order of type $\Delta$ is hereditary, and any hereditary order is isomorphic to a direct product of such orders [CR].

In the rest, assume that $\Lambda$ is an order of type $\Delta$. For any $i\in\Delta_0$, we have a simple object $S_i=S_{i1}\in\fin\Lambda$. For any $j>0$, there exists an indecomposable object $S_{ij}\in\fin\Lambda$ which has the top $S_i$ and the length $j$. It is well-known that $(S_{ij})_{(i,j)\in\Delta_0\times\nnn}$ gives the set of indecomposable objects in $\fin\Lambda$. We denote by $\Pi$ the set of $n$-tuples of partitions. Thus an element $\lambda$ of $\Pi$ can be written as $\lambda=(1^{l_{i1}}2^{l_{i2}}3^{l_{i3}}\cdots)_{i\in \Delta_0}$. Putting $M(\lambda):=\bigoplus_{(i,j)\in \Delta_0\times\nnn}S_{ij}^{l_{ij}}$, we obtain a bijection $M=M_\Lambda:\Pi\rightarrow\fin\Lambda$.

For any $\lambda,\mu,\nu\in\Pi$, there exists $\phi^{\mu}_{\lambda\nu}\in\zzz[T]$ called the {\it Hall polynnomial} such that $F^{M(\mu)}_{M(\lambda)M(\nu)}=\phi^{\mu}_{\lambda\nu}(q_{\Lambda})$ holds for any order $\Lambda$ of type $\Delta$ [Gu]. The {\it generic Hall algebra} $\hall(\Delta)$ is defined as a free $\zzz[T]$-module with the basis $( u_{\lambda})_{\lambda\in\Pi}$ and the multiplication $ u_{\lambda} u_{\nu}=\sum_{\mu\in\Pi}\phi^{\mu}_{\lambda\nu} u_{\mu}$. Then, for any order $\Lambda$ of type $\Delta$, the map $M( u_{\lambda}):= u_{M(\lambda)}$ gives an isomorphism $M:\hall(\Delta)/(T-q_{\Lambda})\rightarrow\hall(\Lambda)$. For the case $n=1$, $\hall(\Delta)$ is the classical Hall algebra [M].

\vskip1em{\bf\ZEBA\ Generic Hall modules }
Let $\Lambda$ be an order of type $\Delta$, $A=\Lambda\otimes_RK$ and $V$ an $A$-module of length $m$. We denote by $\Pi^\infty_m$ the set of $\mu=(m_i)_{i\in\Delta_0}$ with $m_i\in\nnn_{\ge0}$ and $\sum_{i\in\Delta_0}m_i=m$. For any $i\in\Delta_0$, we denote by $P_i$ the projective cover of $S_i$. Putting $M(\mu):=\bigoplus_{i\in\Delta_0}P_i^{m_i}$, we obtain a bijection $M=M_\Lambda:\Pi^\infty_m\rightarrow\overline{\ll}_V(\Lambda)$.

For any $\lambda\in\Pi$ and $\mu,\nu\in\Pi^\infty_m$, there exists $\psi^{\mu}_{\lambda\nu}\in\zzz[T]$ such that $F^{M(\mu)}_{M(\lambda)M(\nu)}=\psi^{\mu}_{\lambda\nu}(q_{\Lambda})$ holds for any order $\Lambda$ of type $\Delta$. An $\hall(\Delta)$-module $\hallm_m(\Delta)$ called the {\it generic Hall module} is defined as a free $\zzz[T]$-module with the basis $( u_{\mu})_{\mu\in\Pi^\infty_m}$ and the action $ u_{\lambda} u_{\nu}=\sum_{\mu\in\Pi^\infty_m}\psi^{\mu}_{\lambda\nu} u_{\mu}$. Then, for any order $\Lambda$ with type $\Delta$ and an $A$-module $V$ of length $m$, the map $M( u_{\mu}):= u_{M(\mu)}$ gives an isomorphism $M:\hallm_m(\Delta)/(T-q_{\Lambda})\rightarrow\hallm_V(\Lambda)$ of $\hall(\Delta)$-modules.

\vskip1em{\bf\ZEBB\ }Put $\hall(\Delta)^{\ccc}_1:=\hall(\Delta)\otimes_{\zzz}\ccc/(T-1)$ and $\hallm_m(\Delta)^{\ccc}_1:=\hallm_m(\Delta)\otimes_{\zzz}\ccc/(T-1)$. Then $\hallm_m(\Delta)^{\ccc}_1$ forms an $\hall(\Delta)^{\ccc}_1$-module. Let $e_i$ be the element of $\dn{gl}_n$ with the $(i,i)$-th entry $1$ and other entries $0$, and $h_i:=e_i-e_{i+1}$. We denote by $\omega_1$ the fundamental weight of $\dn{sl}_n$ defined by $\omega_1(h_i)=\delta_{1i}$ ($1\le i<n$).

\vskip1em{\bf Theorem }{\it
There exists a two-sided ideal $I$ of $\hall(\Delta)^{\ccc}_1$ with the following properties.

(1) $\hall(\Delta)^{\ccc}_1/I$ is isomorphic to the universal enveloping algebra $U(\dn{gl}_n)$ of $\dn{gl}_n$.

(2) For any $m\in\nnn$, the $\hall(\Delta)^{\ccc}_1$-module $\hallm_m(\Delta)^{\ccc}_1$ is annihilated by $I$, and restricts to the ${m+n-1\choose n-1}$-dimensional irreducible $U(\dn{sl}_n)$-module with the highest weight $m\omega_1$.}

\vskip1em
For the case $n=2$, $(\hallm_m(\Delta)^{\ccc}_1)_{m\in\nnn}$ restricts to the set of all finite dimensional irreducible $U(\dn{sl}_2)$-modules. It will be interesting if any irreducible finite dimensional $U(\dn{sl}_n)$-module is constructed by a Hall algebra approach.

\vskip1em{\footnotesize
\begin{center}
{\sc Acknowledgements}
\end{center}

The author would like to thank Professor H. Hijikata for valuable suggestions and continuous encouragement. He also would like to thank Professor K. Nishida, Professor K. W. Roggenkamp and Professor W. Rump for stimulating discussions and encouragement.

\begin{center}
{\bf References}
\end{center}

[A1] M. Auslander: Representation dimension of Artin algebras, Lecture
notes, Queen Mary College, London, 1971.

[A2] M. Auslander: Isolated singularities and existence of almost split sequences. Representation
theory, II (Ottawa, Ont., 1984), 194--242, Lecture Notes in Math., 1178, Springer, Berlin-New York, 1986.

[AB] M. Auslander, R. O. Buchweitz: The homological theory of maximal Cohen-Macaulay approximations, Soc. Math. de France, Mem. no. 38 (1989) 5--37.

[AR] M. Auslander, K. W. Roggenkamp: A characterization of orders of finite lattice type. Invent. Math. 17 (1972), 79--84.

[ARS] M. Auslander, I. Reiten, S. O. Smalo: Representation theory of Artin algebras. Cambridge Studies in Advanced Mathematics, 36. Cambridge University Press, Cambridge, 1995. 

[AS1] M. Auslander, S. O. Smalo: Preprojective modules over Artin algebras. J. Algebra 66 (1980), no. 1, 61--122. 

[AS2] M. Auslander, S. O. Smalo: Almost split sequences in subcategories. J. Algebra 69 (1981), no. 2, 426--454.

[As] H. Asashiba: Realization of general and special linear algebras via Hall algebras, preprint.

[B] H. Bass: Finitistic dimension and a homological generalization of semiprimary rings, Trans. Amer. Math. Soc. 95 (1960) 466-488.

[BR1] C. J. Bushnell, I. Reiner: Zeta functions of arithmetic orders and Solomon's conjectures. Math. Z. 173 (1980), no. 2, 135--161.

[BR2] C. J. Bushnell, I. Reiner: Zeta-functions of orders. Integral representations and applications (Oberwolfach, 1980), pp. 159--173, Lecture Notes in Math., 882, Springer, Berlin-New York, 1981.

[BR3] C. J. Bushnell, I. Reiner: The prime ideal theorem in noncommutative arithmetic. Math. Z. 181 (1982), no. 2, 143--170.

[BR4]  C. J. Bushnell, I. Reiner: New asymptotic formulas for the distribution of left ideals of orders. J. Reine Angew. Math. 364 (1986), 149--170.

[CPS1] E. Cline, B. Parshall, L. Scott: Finite-dimensional algebras and highest weight categories. J. Reine Angew. Math. 391 (1988), 85--99.

[CPS2] E. Cline, B. Parshall, L. Scott: Algebraic stratification in representation categories. J. Algebra 117 (1988), no. 2, 504--521.

[CR] C. W. Curtis, I. Reiner: Methods of representation theory. Vol. I. With applications to finite groups and orders. A Wiley-Interscience Publication. John Wiley \& Sons, Inc., New York, 1990.

[D] M. Denert: Solomon's second conjecture: a proof for local hereditary orders in central simple algebras. J. Algebra 139 (1991), no. 1, 70--89.

[DK1] Y. A. Drozd, V. V. Kiri\v cenko: The quasi-Bass orders. (Russian) Izv. Akad. Nauk SSSR Ser.
Mat. 36 (1972), 328--370.

[DK2] Y. A. Drozd, V. V. Kiri\v cenko:
Primary orders with a finite number of indecomposable
representations. (Russian) Izv. Akad. Nauk SSSR Ser. Mat. 37 (1973), 715--736.

[DKR] Y. A. Drozd, V. V. Kiri\v cenko, A. V. Ro\u\i ter:
Hereditary and Bass orders. (Russian) Izv.
Akad. Nauk SSSR Ser. Mat. 31 1967 1415--1436.

[DR1] V. Dlab, C. M. Ringel: Quasi-hereditary algebras. 
Illinois J. Math. 33 (1989), no. 2, 280--291. 

[DR2] V. Dlab, C. M. Ringel: Every semiprimary ring is the endomorphism ring of a projective module over a quasihereditary ring. Proc. Amer. Math. Soc.
107 (1989), no. 1, 1--5.

[DR3] V. Dlab, C. M. Ringel: Auslander algebras as quasi-hereditary algebras. J. London Math. Soc. (2) 39 (1989), no. 3, 457--466.

[DW] E. Dieterich, A. Wiedemann: The Auslander-Reiten quiver of a simple curve singularity, Trans. Amer. Math. Soc. 294 (1986) 455-475.

[EHIS] K. Erdmann, T. Holm, O. Iyama, J. Schr\"oer: Radical embeddings and representation dimension, to appear in Advances in Mathematics.

[FGR] R. M. Fossum, P. Griffith, I. Reiten: Trivial extensions of abelian categories. Lecture Notes in Mathematics, Vol. 456. Springer-Verlag, Berlin-New York, 1975.

[G] V. M. Galkin: Zeta-functions of certain one-dimensional rings. (Russian) Izv. Akad. Nauk SSSR Ser. Mat. 37 (1973), 3--19.

[Gu] J. Guo: The Hall polynomials of a cyclic serial algebra. Comm. Algebra 23 (1995), no. 2, 743--751.

[HN1] H. Hijikata, K. Nishida:
Classification of Bass orders. J. Reine Angew. Math. 431 (1992),
191--220.

[HN2] H. Hijikata, K. Nishida:
Bass orders in nonsemisimple algebras. J. Math. Kyoto Univ. 34 (1994),
no. 4, 797--837.

[HN3] H. Hijikata, K. Nishida:
Primary orders of finite representation type. J. Algebra 192
(1997), no. 2, 592--640. 

[H] T. Holm: Representation dimension of some tame blocks of finite groups, to appear in Algebra Colloquium.

[I1] O. Iyama: A generalization of rejection lemma of Drozd-Kirichenko. J. Math. Soc. Japan 50 (1998), no. 3, 697--718.

[I2] O. Iyama: Some categories of lattices associated to
a central idempotent. J. Math. Kyoto Univ. 38 (1998), no. 3, 487--501.

[I3] O. Iyama: $\tau$-categories I--III, to appear in Algebras and Representation theory.

[I4] O. Iyama: A proof of Solomon's second conjecture on local zeta functions of orders. J. Algebra 259 (2003), no. 1, 119--126.

[I5] O. Iyama: Finiteness of Representation dimension. Proc. Amer. Math. Soc. 131 (2003), 1011-1014.

[I6] O. Iyama: Symmetry and duality on $n$-Gorenstein rings, to appear in Journal of Algebra.

[I7] O. Iyama: On Hall algebras of hereditary orders, to appear in Communications in Algebra.

[I8] O. Iyama: The relationship between homological properties and representation theoretic realization of artin algebras, preprint.

[I9] O. Iyama: On the representation dimension of orders, in preparation.

[IT] K. Igusa, G. Todorov: On the finitistic global dimension conjecture, preprint.

[KW] S. K\"onig, A. Wiedemann: Global dimension two orders are quasi-hereditary. Manuscripta Math. 66 (1989), no. 1, 17--23.

[L] S. Lang: Algebraic number theory. Reading, Massachusetts: Addison-Wesley 1970.

[Lu] G. Lusztig: Quivers, perverse sheaves, and quantized enveloping algebras, J. Amer. Math. Soc., 4 (1991), 365-421.

[M] I. G. Macdonald: Symmetric functions and Hall polynomials. Second edition. With contributions by A. Zelevinsky. Oxford Mathematical Monographs. Oxford Science Publications. The Clarendon Press, Oxford University Press, New York, 1995. x+475 pp.

[R1] C. M. Ringel: Hall algebras. Topics in algebra, Part 1 (Warsaw, 1988), 433--447, Banach Center Publ., 26, Part 1, PWN, Warsaw, 1990.

[R2] C. M. Ringel: The composition algebra of a cyclic quiver. Proc. London Math. Soc. (3) 66 (1993), no. 3, 507--537.

[Ro1] K. W. Roggenkamp: Lattices over orders. II. Lecture Notes in Mathematics, Vol. 142 Springer-Verlag, Berlin-New York 1970.

[Ro2] K. W. Roggenkamp: The construction of almost split sequences for integral group rings and orders. Comm. Algebra 5 (1977), 1363-1373.

[Ru1] W. Rump: Ladder functors with an application to representation-finite Artinian rings. To Mirela \c Stef\u anescu, at her 60's. An. \c Stiin\c t. Univ. Ovidius Constan\c ta Ser. Mat. 9 (2001), no. 1, 107--123. 

[Ru2] W. Rump: Lattice-finite rings, preprint.

[Ru3] W. Rump: The category of lattices over a lattice-finite ring, preprint.

[S1] L. Solomon: Zeta functions and integral representation theory. Advances in Math. 26 (1977), no. 3, 306--326.

[S2] L. Solomon: Partially ordered sets with colors. Relations between combinatorics and other parts of mathematics (Proc. Sympos. Pure Math., Ohio State Univ.,
Columbus, Ohio, 1978), pp. 309--329, Proc. Sympos. Pure Math., XXXIV, Amer. Math. Soc., Providence, R.I., 1979.

[Sc] O. Schiffmann: The Hall algebra of a cyclic quiver and canonical bases of Fock spaces. Internat. Math. Res. Notices 2000, no. 8, 413--440.

[T] H. Tachikawa: Quasi-Frobenius rings and generalizations. Lecture Notes in Mathematics, Vol. 351. Springer-Verlag, Berlin-New York, 1973.

[W] A. Wiedemann: Classification of the Auslander-Reiten quivers of local Gorenstein orders and a characterization of the simple curve singularities. J. Pure Appl. Algebra 41 (1986), no. 2-3, 305--329.

[X1] C. C. Xi: On the representation dimension of finite dimensional algebras. J. Algebra 226 (2000), no. 1, 332--346.

[X2] C. C. Xi: Representation dimension and quasi-hereditary algebras. Adv. Math. 168 (2002), no. 2, 193--212.

[X3] C. C. Xi: The relative transpose of a module, preprint.

[X4] C. C. Xi: The finitistic dimension conjecture and representation-finite algebras, preprint.

[Y] Y. Yoshino: Cohen-Macaulay modules over Cohen-Macaulay rings. London Mathematical Society Lecture Note Series, 146. Cambridge University Press, Cambridge, 1990. viii+177 pp.

[Z] H. B. Zimmermann: The finitistic dimension conjectures---a tale of $3.5$ decades. Abelian groups and modules (Padova, 1994), 501--517, Math. Appl., 343, Kluwer Acad. Publ., Dordrecht, 1995.
\end{document}